\documentclass[article]{article}
\usepackage{color}
\usepackage{graphicx}
\usepackage{amsmath}
\usepackage{amssymb}
\usepackage{mathrsfs}
\usepackage[numbers,sort&compress]{natbib}
\usepackage{amsthm}
\usepackage{tikz}
\usepackage{float}
\newtheorem{theorem}{Theorem}

\newtheorem{remark}{Remark}

\usepackage {caption}
\usepackage{float}
\usepackage{subfigure}
\title{Critical edge behavior  in the singularly perturbed Pollaczek-Jacobi type unitary ensemble  }
\date{}
\author{Zhaoyu Wang,  Engui Fan$^1$\thanks{\ Corresponding author and email address: faneg@fudan.edu.cn } }
\footnotetext[1]{ \  School of Mathematical Sciences, Fudan University, Shanghai 200433, P.R. China.}
\begin{document}
\maketitle
\textbf{Abstract:} In this paper, we study the  strong asymptotic  for   the  orthogonal polynomials and  universality
associated with singularly perturbed  Pollaczek-Jacobi type   weight
$$w_{p_J2}(x,t)=e^{-\frac{t}{x(1-x)}}x^\alpha(1-x)^\beta,  $$
  where $t \ge 0$, $\alpha >0$, $\beta >0$ and $x \in [0,1].$   Our  main results obtained here  include two aspects:
  { I. Strong asymptotics:}   We obtain the strong asymptotic expansions for the monic Pollaczek-Jacobi type
  orthogonal polynomials in different interval $(0,1)$ and outside of  interval  $\mathbb{C}\backslash (0,1)$, respectively;
Due to the effect of   $\frac{t}{x(1-x)}$ for varying $t$,  different asymptotic behaviors   at the hard edge $0$ and $1$
 were found with different scaling schemes.  Specifically, the uniform asymptotic behavior can be expressed as a
 Airy function in the neighborhood of point $1$ as $\zeta= 2n^2t \to \infty, n\to \infty$,  while it is given by a Bessel function as $\zeta \to 0, n \to \infty$.
 { II. Universality:}  We respectively  calculate   the limit of  the eigenvalue correlation kernel    in the bulk of the spectrum and  at the both side of hard edge,
  which will  involve  a $\psi$-functions  associated with a particular Painlev$\acute{e}$ \uppercase\expandafter{\romannumeral3} equation near $x=\pm 1$.
    Further, we also prove the $\psi$-funcation can be approximated by a Bessel kernel as $\zeta \to 0$ compared with a Airy kernel as $\zeta \to \infty$.
    Our analysis   is  based on the Deift-Zhou nonlinear steepest descent method for  the  Riemann-Hilbert problems.
\\ \textbf{Kerwords:}  singularly perturbed  Pollaczek-Jacobi type  weight,  strong asymptotic expansions, eigenvalue correlation  kernel,  universality,  Riemann-Hilbert problem, Deift-Zhou nonlinear steepest descent method.

\newpage
\numberwithin{equation}{section}
\section{Introduction }

   The orthogonal polynomials (OPS)  associated with  different  weights  have  important applications  in random matrix theory \cite{Deift1999Orthogonal,mehta2004random},  approximation theory \cite{Szeg1994Orthogonal} and the zeros distribution theory \cite{Chihara1978An}.
    There are many ways to classify the weight functions.  One way, just like  the classification of random matrices ensembles, is to identify the weights into
      Hermite weight, Gaussian weight, Laguerre weight and  Jacobi weight and so on \cite{Forrester2010Log}.  Another way is to classify  the weights into Szeg$\ddot{o}$  class   or non-Szeg$\ddot{o}$
      class according the following Szeg$\ddot{o}$ condition \cite{Szeg1994Orthogonal}.
    A weight $w(x)$ is said to be of the Szeg$\ddot{o}$ class if   it   has a   definition in the domain $[-1,1]$   and   satisfies
    $$\int_{-1}^{1} \frac{\text{ln}w(x)}{\sqrt{1-x^2}} \,\mathrm{d} x > -\infty.$$
  It is interest to study the weight in the Szeg$\ddot{o}$ class.  It has been shown  that  a large number  of weights belong to the Szeg$\ddot{o}$ class.
 Especially, the Jacobi weight is  a typical example of the Szeg$\ddot{o}$ class, in details  see \cite{kuijlaars2001riemann}.

  The  asymptotics of OPS  with  respect to  the classical Jacobi weight has been studied in \cite{Szeg1994Orthogonal}.
    There are different  generalizations to  the classical  Jacobi weight.  For example, Kuijlaars et al  investigated  the OPS
    about the modified Jacobi weight  \cite{kuijlaars2001riemann}
    $$w(x)=(1-x)^\alpha(1+x)^\beta h(x), \quad x \in (-1,1),$$
    where $\alpha, \beta>-1$ and   $h(x)>0$  is a  real analytic function  for $x \in (-1,1)$.   They further obtained   the  asymptotic behavior for the Hankel determinants,
    the recurrence coefficients the leading coefficients of the OPS.
     Xu and Zhao have obtained the asymptotics of the polynomials orthogonal with respect to the perturbed Chebyshev weight in \cite{xu2013critical},
    $$w(x)=(1-x)^{-\frac{1}{2}}(t_n^2-x^2)^\alpha, \quad x \in (-1,1), t_n>1.$$
    At the same time, they give the critical edge behavior for this weight with different Painlev$\acute{e}$ equation according to the different property of $t_n$.

      However, there are some weights that furnish a non-Szeg$\ddot{o}$ class such as  the Pollaczek weight which will considered in this paper.
       The Pollaczek polynomials are orthogonal to the following Pollaczek weight function
    $$\begin{aligned}
        w(x,a,b) & =\text{exp}(2\theta-\pi)h(\theta)[\text{cosh}(\pi h(\theta))^{-1}],
     \end{aligned}$$
    where $\theta=\text{arccos}\sqrt{x}$, $h(\theta)=(a \text{cos}\theta+b)/(2\text{sin}\theta)$, and a and b are  constants, and $|b|<a$  \cite{Szeg1994Orthogonal,ZhouUniform}.
   Compared by the Jacobi polynomials, the Pollaczek polynomials show a singular behavior in some items, such as the properties of the zero and the gap   \cite{Szeg1994Orthogonal}. Ismail have studied the large-n asymptotic behavior of the Pollaczek polynomials using the confluent Horn function in \cite{ismail1994asymptotics}.    Bo and Wong have investigated the asymptotics with the technique of the Airy function, see \cite{rui2016asymptotic}.

    \par In recent years,   there are some interests  on the study on the associated weight functions with singular behavior in several different areas of mathematics and physics. For instance, Berry and Shukla  studied the statistics for the zeros of the Riemann Zeta function in the Gaussian unitary ensemble    of random matrices, with the singularly perturbed Gaussion weight \cite{BerryTuck}
    $$w(x;z,s)=\text{exp}(-\frac{z^2}{2x^2}+\frac{s}{x}-\frac{x^2}{2}), z \in \mathbb{R}\setminus {0}, s\in [0,\infty), x \in \mathbb{R}.$$
     Then, Brightmore, Mezzadri and Mo studied the double scaling limit as $N \to \infty$ and computed the asymptotics of the partition function when $z$ and $t$ are of order $\mathcal{O}(N^{-1/2})$  which can be characterized   as a  particular Painlev$\acute{e}$ \uppercase\expandafter{\romannumeral3} equation \cite{BrightmoreA}.
     Besides, Lukyanov  calculated the finite temperature expectation values of the exponential fields to one-loop order of the semi-classical expansion under the singularly perturbed Laguerre weight
     $$w(x,t)=x^\alpha e^{-x-t/x}, \alpha>0, t>0, x\in (0,\infty),$$
     which is the Laguerre weight $x^\alpha e^{-x}$ perturbed by a multiplicative factor $e^{s/x}$ inducing an infinitely strong zero at the origin \cite{Lukyanov2000Finite}. In 2010, Chen and Its studied a certain linear statistics of the unitary Laguerre ensembles and the determination of the associated Hankel determinant and recurrence coefficients \cite{YangPainlev}.
     After that, Xu, Dai, and Zhao showed the eigenvalue correlation kernel associated with a particular Painlev$\acute{e}$ \uppercase\expandafter{\romannumeral3} equation.
      They also studied the transition of this limiting kernel to the Bessel and Airy kernels \cite{XuCritical}.
    Further, Chen and Dai introduced a sequence of polynomials orthogonal with respect to a one-parameter family of   the Pollaczek-Jacobi type weight
   \begin{align}
    w(x,t)=x^{\alpha}(1-x)^{\beta} e^{-t/x}, \quad t \ge 0, x \in [0,1].\label{1singl}
    \end{align}
     If $t=0$, this becomes to a shifted Jacobi weight. They expressed the recurrence coefficients in terms of a set of auxiliary quantities which equaled to a particular Painlev$\acute{e}$ \uppercase\expandafter{\romannumeral5} or allied functions using the ladder operator \cite{Chen2010Painlev}. Chen and coauthhors showed the uniform asymptotic expansions for the moninc OPS and the recurrence coefficients and leading coefficient of the OPS were expressed in terms of a particular Painlev$\acute{e}$ \uppercase\expandafter{\romannumeral3} using the Riemann-Hilbert approach \cite{ChenThe}.  It is easy to check  that this weight belongs to non-Szeg$\ddot{o}$ since this weight is in some extent more ``singular".

    \par The above weights  considered admit   the first order singularities.  However, it is also interesting to study  weights   with higher order singularities.
    Atkin, Claeys and Mezzadri  obtained asymptotics for the partition functions associated to the Laguerre and Gaussian Unitary Ensembles perturbed with a pole of order $k$ at the origin and their results were described in terms of a hierarchy of higher order analogs to the PIII equation \cite{AtkinRandom}.
     Also, Dan, Xu and Zhang  introduced the Fredholm determinant of an integrable operator whose kernel is constructed out of the $\Psi$-function associated with a hierarchy of higher order analogues to the Painlev$\acute{e}$ \uppercase\expandafter{\romannumeral3}equation and obtained the large $s$ asymptotics of the Fredholm determinant using the RH method \cite{Dai2017Gap}.

    \par One of  interesting phenomenon in random matrices is that as  the order of the matrix goes to be  large,
    the eigenvalues of different types of random matrix ensembles show the same statistical law. This phenomenon is known as universality, see \cite{deift1999uniform,magnus1986freud}.
   For example, for  the Jacobi unitary ensemble,  its   limiting mean eigenvalue density admits the property
     $${\lim_{n \to +\infty}}\frac{1}{n}K_n(x,x)=\frac{1}{\pi\sqrt{1-x^2}},\ \ \ \  x \in (-1,1).$$
    The  limit  of  $K_n$ is given by the sine kernels in the bulk of the spectrum \cite{deift1999uniform,lubinsky2008universality},
    \begin{equation}\label{sins}
        \begin{aligned}
            \mathbb{S}(x,y)=\frac{\text{sin}\pi(x-y)}{x-y}.
        \end{aligned}
    \end{equation}
  And  the  limit  of  $K_n$  shows   the Bessel kernel at the hard edge of the spectrum \cite{mehta2004random}
    \begin{equation}\label{besselk}
        \begin{aligned}
            \mathbb{J}_\alpha(x,y)=\frac{J_\alpha(\sqrt{x})\sqrt{y}J_\alpha'(\sqrt{y})-J_\alpha(\sqrt{y})\sqrt{x}J_\alpha'(\sqrt{x})}{2(x-y)}.
        \end{aligned}
    \end{equation}

    \par In this paper, we apply  the Deift-Zhou nonlinear steepest method to  study  the   OPS associated with the Pollaczek-Jacobi type
    weight involving  two singularities at the edge
    \begin{equation}\label{wpj2}
       \begin{aligned}
         w_{p_J2}(x) = x^{\alpha}(1-x)^{\beta} e^{-t/x(1-x)},
       \end{aligned}
    \end{equation}
   where $  t \ge 0, \   \alpha >0,\   \beta >0,\   x \in  \left[ 0, 1  \right]$.
  Compared  to the weight (\ref{1singl}),   our  weight (\ref{wpj2})  involves  one  more  singular  perturbation  $1/(z-1)$.
  Our  purpose is to see   whether  some  new  phenomena  occur     with  the  change of  singularities.
 For our  weight (\ref{wpj2}),  we have to  construct another model RH problem to match the original  RH problem,
     rather than simply using the RH problem for the Bessel model in \cite{ChenThe}.
    To build a suitable model, using the same method in \cite{XuCritical,xu2015critical},  we formulate a model RH problem, which is equivalent to a Cauchy-type singular integral equations, see \cite{fokas2006painleve}.  Moreover, there is an interesting result that this model can transform from a modified Bessel model RH problem as $s \to 0,\, n\to \infty$ to another modified Airy function as $s \to \infty,\, n\to \infty$.

    This paper is organized as follows. In Section 2, we state our main results which will be
    presented in this paper.  In section 3, we propose a RH problem with respect to the
    weight (\ref{wpj2})   to characterize   the   Pollaczek-Jacobi type  OPS.
     In Section 3, we apply the Deift-Zhou nonlinear steepest descent method to analyze
the RH problem.  In section 4, based on the asymptotic analysis, we   provide    detail   proofs to all results stated in the section 2.

 \section{ Statement of Main Results}

     \par Let introduce some notations used in this paper. Suppose that  $\pi_n(x)=\pi_n(x;w_{p_J2})$ denote the monic  OPS
    of degree $n$ orthogonal with respect to the weight (\ref{wpj2}), then it holds  that
    \begin{equation}
       \begin{aligned}
        \int_{0}^{1} \pi_n(x) \pi_m(x) w_{p_J2}(x,t) \,\mathrm{d} x =h_n(t)\delta_{nm}.
       \end{aligned}
    \end{equation}
     Let $P_n(x)=P_n(x;w_{p_J2}(x,t))$  to  denote  normalized  polynomials with respect to the weight $w_{p_J2}(x,t)$, then we have
     $$P_n(x)=\gamma_n\pi_n(x),$$
     where $\gamma_n>0$ is the leading coefficient of $P_n(x)$ and $h_n=\gamma_n^{-2}$.
    \par According to basic property of orthogonal  polynomials \cite{Deift1999Orthogonal,Forrester2010Log,mehta2004random}, then  $\pi_n(x)$ meets the three-term recurrence relation
    \begin{equation}
        \begin{aligned}
          x\pi_n(x)=\pi_{n+1}(x)+\alpha_n(t)\pi_n(x)+\beta_n(t)\pi_{n-1}(x),
        \end{aligned}
    \end{equation}
    and the eigenvalue correlation kernel is expressed as
    \begin{equation}
        \begin{aligned}
           K_n(x,y;t)=\sqrt{w_{p_J2}(x,t)}\sqrt{w_{p_J2}(y,t)} \sum_{i=1}^{n}P_i(x)P_i(y).
        \end{aligned}
    \end{equation}
    Using the Christoffel-Darboux formula  \cite{Szeg1994Orthogonal}, we can get
    \begin{equation}\label{kn}
        \begin{aligned}
            K_n(x,y;t)=\gamma_{n-1}^2\sqrt{w_{p_J2}(x,t)}\sqrt{w_{p_J2}(y,t)} \frac{\pi_n(x)\pi_{n-1}(y)-\pi_n(y)\pi_{n-1}(x)}{(x-y)}.
        \end{aligned}
    \end{equation}

  \subsection{Asymptotic for monic OPS }
     \par   Using the Deift-Zhou steepest descent method, we get the uniform asymptotic expansion for the monic Pollaczek-Jacobi type OPS.

     \begin{theorem}
     For  $z\in  \mathbb{C} \setminus \left[ 0, 1  \right]$,  as $n \to \infty$,   $\pi_n(z)$ has the following asymptotic approximation
    \begin{small}
      \begin{equation}
        \begin{aligned}
            \pi_n(z)=2^{-(2n+\alpha+\beta+1)}(\varphi(z))^{\frac{2n+\alpha+\beta+1}{2}}(z-1)^{-\frac{2\beta+1}{4}}z^{-\frac{2\alpha+1}{4}}e^{\frac{t}{2z(1-z)}}(1+\mathcal{O}(\frac{1}{n})),
        \end{aligned}
      \end{equation}
    \end{small}
    where $t>0$, $\varphi(z)=2z-1+2\sqrt{z(z-1)}$ and the error terms have a uniform for z on compact subsets of $\bar{\mathbb{C}}\setminus [0, 1]$ as well as a full asymptotic expansion in terms of $1/n$.
\end{theorem}

    \par We also calculate the strong asymptotic expansion in the bulk of the orthogonality interval $(0, 1)$.

    \begin{theorem}
     Let $\mathcal{K}_1$ be a compact subset of $(0,1)$. When $0<x<1$ and $n \to \infty$,
     \begin{small}
      \begin{equation}
       \begin{aligned}
           \pi_n(x)= \frac{2^{-(\alpha+\beta+2n)}}{\sqrt{w_{p_J2}(x,t)}x^{\frac{1}{4}}(1-x)^{\frac{1}{4}}}\left[{ (1+\mathcal{O}(\frac{1}{n}))\text{cos}((2n+\alpha+\beta+1)\text{arccos}\sqrt{x}-\frac{\beta \pi }{2}-\frac{\pi}{4})} \right.\\
           \left. {+\mathcal{O}(\frac{1}{n})\text{cos}((2n+\alpha+\beta+1)\text{arccos}\sqrt{x}-\frac{\beta \pi }{2}+\frac{\pi}{4})} \right]
       \end{aligned}
      \end{equation}
     \end{small}
    The above expansion is uniform for $x \in \mathcal{K}_1$.
\end{theorem}

    \par As $\zeta=2n^2t \to 0$, the asymptotic behavior of $\pi_n$ near the endpoint 1 is described by the following theorem.

 \begin{theorem}
     There is a small and positive constant $r$  which satisfies $0< r \ll \frac{1}{2}$,  such that $x \in (1-r,1)$,  as  $n\to \infty$, then
      \begin{small}
       \begin{equation}
            \begin{aligned}
                \pi_n(x) &= \frac{\sqrt{\pi}(2n\text{arccos}\sqrt{x})^{\frac{1}{2}}}{2^{2n+\alpha+\beta}x^{\frac{1}{4}}(1-x)^{\frac{1}{4}}\sqrt{2w_{p_J2}(x,t)}} e^{\frac{-\zeta}{(2n\text{arccos}\sqrt{x})^2}}\\
                         & \times [(1+\mathcal{O}(1/n)+\mathcal{O}(\zeta))(\text{cos}\theta_1 J_\beta (2n\text{arccos}\sqrt{x})+\text{sin}\theta_1 J'_\beta (2n\text{arccos}\sqrt{z}))\\
                         & + (\mathcal{O}(1/n)+\mathcal{O}(\zeta))(\text{cos}\theta_2 J_\beta (2n\text{arccos}\sqrt{x})+\text{sin}\theta_2 J'_\beta (2n\text{arccos}\sqrt{z}))],
            \end{aligned}
        \end{equation}
      \end{small}
        where
        \begin{equation} \label{t1t2}
            \begin{aligned}
                \theta_1=(\alpha+\beta+1)\text{arccos}\sqrt{x},\; \text{and}\; \theta_2=(\alpha+\beta-1)\text{arccos}\sqrt{x},
            \end{aligned}
        \end{equation}
        and $J_\beta$ is the Bessel function of order $\beta$. The error terms hold uniformly for $x$ in compact subsets of $(1-r,1)\cap \mathbb{D}_1$,
         $\mathbb{D}_1=\{x|2n\pi \text{arccos}\sqrt{x} \in (0,\infty), n\to \infty \}$.
\end{theorem}

    \par Near the endpoint 1, the strong asymptotic expansion can be calculated in terms of the Airy function as $\zeta=2n^2t\to \infty$.

       \begin{theorem}
   There is a small and fixed constant $\eta$, $0<\eta \ll \frac{1}{2}$, such that $x \in (1-\eta, 1)$ and $\zeta= 2n^2t \to \infty$,
      \begin{small}
       \begin{equation}
        \begin{aligned}
            \pi_n(x) &=\frac{\sqrt{\pi} (2n\text{arccos}\sqrt{x})^{\frac{1}{2}}}{2^{2n+\alpha+\beta}x^{\frac{1}{4}}(1-x)^{\frac{1}{4}}\sqrt{w_{p_J2}(x,t)}}\\
                     & \times \{ (1+\mathcal{O}(n^{-\frac{2}{3}})) [(\text{cos}\theta_1 d_{11}+ \text{sin}\theta_1 (d_{21}-b d_{11})(2n\text{arccos}\sqrt{x})^{-1})Ai(\tilde{\xi})   \\
                     & +(\text{cos}\theta_1 d_{12}- \text{sin}\theta_1 (d_{22}+b d_{12})(2n\text{arccos}\sqrt{x})^{-1})Ai'(\tilde{\xi})]  \}\\
                     & + \mathcal{O}(n^{-\frac{2}{3}})[(\text{cos}\theta_2 d_{11}+ \text{sin}\theta_2 (d_{21}-b d_{11})(2n\text{arccos}\sqrt{x})^{-1})Ai(\tilde{\xi})   \\
                     & +(\text{cos}\theta_2 d_{12}- \text{sin}\theta_2 (d_{22}+b d_{12})(2n\text{arccos}\sqrt{x})^{-1})Ai'(\tilde{\xi})]  \},
        \end{aligned}
       \end{equation}
       \end{small}
       where $\theta_1, \theta_2$ are in (\ref{t1t2}), $\tilde{\xi}=(\frac{3}{2})^{\frac{2}{3}}(1-|\lambda|)\zeta^{\frac{2}{9}}|\lambda|^{-\frac{2}{3}}$,
       $\lambda=\zeta^{-\frac{2}{3}}n^2f_1(z)$, $f_1(z)=(\text{log}\varepsilon(z))^2$.  The error terms hold uniformly for $x$
        in compact subsets of $(1-\eta,1)\cap \mathbb{D}_2,$
        where
       $$\mathbb{D}_2=\{ x| (\frac{3}{2})^{\frac{2}{3}} (1-|\lambda|) \zeta^{\frac{2}{9}} \in (-\infty, 0), |\lambda|=\zeta^{-\frac{2}{3}}
       (2n\text{arccos}\sqrt{x})^2, \zeta \to \infty, n\to \infty \}.$$
    \end{theorem}

       At the same time, $b, d_{11}, d_{12}, d_{21}, \text{and}\; d_{22}$ are defined as follows:
       \begin{equation}\label{b}
           \begin{aligned}
               b = \frac{7}{72(|\lambda|-1)}
           \end{aligned}
       \end{equation}
       \begin{equation}\label{d11}
           \begin{aligned}
               d_{11}= (\frac{3}{2})^{\frac{1}{6}} \zeta^{-\frac{1}{9}}|\lambda|^{-\frac{1}{6}} \text{cos}(\frac{\beta}{2}\text{arccos}\frac{1}{\sqrt{|\lambda|}}),
             \end{aligned}
        \end{equation}
        \begin{equation}\label{d22}
           \begin{aligned}
               d_{22}= (\frac{3}{2})^{-\frac{1}{6}} \zeta^{\frac{1}{9}}|\lambda|^{\frac{1}{6}} \text{cos}(\frac{\beta}{2}\text{arccos}\frac{1}{\sqrt{|\lambda|}}),
             \end{aligned}
        \end{equation}
        \begin{equation}\label{d12}
           \begin{aligned}
               d_{12}= (\frac{3}{2})^{-\frac{1}{6}} \zeta^{-\frac{2}{9}}(|\lambda|-1)^{-\frac{1}{2}}|\lambda|^{\frac{1}{6}} \text{sin}(\frac{\beta}{2}\text{arccos}\frac{1}{\sqrt{|\lambda|}}),
             \end{aligned}
        \end{equation}
        \begin{equation}\label{d21}
           \begin{aligned}
               d_{21}= (\frac{3}{2})^{\frac{1}{6}} \zeta^{\frac{2}{9}}(|\lambda|-1)^{\frac{1}{2}}|\lambda|^{-\frac{1}{6}} \text{sin}(\frac{\beta}{2}\text{arccos}\frac{1}{\sqrt{|\lambda|}}).
             \end{aligned}
        \end{equation}

    \par The strong asymptotic expansion is expressed by the Bessel function of order $\alpha$ as $\zeta = 2n^2 t \to 0$ in the neighborhood of the point 0.

    \begin{theorem}
    There exists a small and positive constant $\delta$, $0 < \delta \ll \frac{1}{2} $ so that for $x \in (0,\delta)$ and $\zeta = 2 n^2 t$, when $n \to \infty$, $\zeta \to 0$,
    \begin{small}
     \begin{equation}\nonumber
        \begin{aligned}
            \pi_n(x) &= \frac{\sqrt{\pi}(n(\pi - 2\text{arccos}\sqrt{x}))^{\frac{1}{2}}}{2^{2n+\alpha+\beta}x^{\frac{1}{4}}(1-x)^{\frac{1}{4}}\sqrt{2w_{p_J2}(x,t)}} e^{\frac{-\bar{\zeta}}{(n(\pi - 2\text{arccos}\sqrt{x}))^2}}\\
                     & \times [(1+\mathcal{O}(1/n)+\mathcal{O}(\zeta))(\text{sin}\theta_3 J_\alpha (n(\pi - 2\text{arccos}\sqrt{x}))+\text{cos}\theta_3 J'_\alpha (n(\pi - 2\text{arccos}\sqrt{z})))\\
                     & + (\mathcal{O}(1/n)+\mathcal{O}(\zeta))(\text{sin}\theta_4 J_\alpha (n(\pi - 2\text{arccos}\sqrt{x}))+\text{cos}\theta_4 J'_\alpha (n(\pi - 2\text{arccos}\sqrt{z})))],
        \end{aligned}
    \end{equation}
   \end{small}
    where
    \begin{equation}\label{t3}
        \begin{aligned}
           \theta_3= (\alpha+\beta+1)\text{arccos}\sqrt{x}-\frac{(\alpha+\beta)\pi}{2}
        \end{aligned}
     \end{equation}
     \begin{equation}\label{t4}
        \begin{aligned}
           \theta_4= (\alpha+\beta-1)\text{arccos}\sqrt{x}-\frac{(\alpha+\beta)\pi}{2}
        \end{aligned}
     \end{equation}
     and $J_\alpha$ is the Bessel function of order $\alpha$.
    The error terms is uniformly for $x$ in compact subsets of $(0,\delta)\cap \mathbb{D}_3$, where $\mathbb{D}_3=\{ x| n(\pi - 2\text{arccos}\sqrt{x})\in (0,\infty), n\to \infty \}$.
\end{theorem}

    \par We can obtain similar asymptotic behavior in terms of Airy function of order $\alpha$ near the endpoint 0 as $\zeta= 2n^2 t \to \infty$.

\begin{theorem} There exists a small and fixed $\epsilon$, $0< \epsilon \ll \frac{1}{2} $, so that $x \in (0, \epsilon)$, if $n\to \infty,\; \zeta= 2n^2 t \to \infty$,
    \begin{small}
     \begin{equation}
        \begin{aligned}
            \pi_n(x) &= \frac{\sqrt{\pi}(n(\pi - 2\text{arccos}\sqrt{x}))^{\frac{1}{2}}}{2^{2n+\alpha+\beta}x^{\frac{1}{4}}(1-x)^{\frac{1}{4}}\sqrt{2w_{p_J2}(x,t)}} \\
                     & \times \{ [1+\mathcal{O}(n^{-\frac{2}{3}})] [(\text{sin}\theta_3d_{11}+ \text{cos}\theta_3(n(\pi - 2\text{arccos}\sqrt{x}))^{-1}(bd_{11}+d_{21}))Ai(\hat{\xi})\\
                     & + (\text{sin}\theta_3d_{11}+ \text{cos}\theta_3(n(\pi - 2\text{arccos}\sqrt{x}))^{-1}(bd_{12}-d_{22}))Ai'(\hat{\xi})]\\
                     & + \mathcal{O}(n^{-\frac{2}{3}})[(\text{sin}\theta_4d_{11}+ \text{cos}\theta_4(n(\pi - 2\text{arccos}\sqrt{x}))^{-1}(bd_{11}+d_{21}))Ai(\hat{\xi})\\
                     & + (\text{sin}\theta_4d_{11}+ \text{cos}\theta_4(n(\pi - 2\text{arccos}\sqrt{x}))^{-1}(bd_{12}-d_{22}))Ai'(\hat{\xi})]\},
        \end{aligned}
     \end{equation}
    \end{small}
    where $Ai$ is the Airy function, $\hat{\xi}=(\frac{3}{2})^{\frac{2}{3}}(1-|\bar{\lambda}|)\zeta^{\frac{2}{9}}|\bar{\lambda}|^{-\frac{2}{3}}$, $\bar{\lambda} = \zeta^{-\frac{2}{3}}n^2 f_0(z)$ with $f_0(z)=(\text{log}\frac{\varphi(z)}{\varphi(0)})^2$. The error terms hold uniformly for $x \in \mathcal{K}_2$, where $\mathcal{K}_2$ is any compact subset of $(0, \epsilon) \cap \mathbb{D}_4, t \in (0, c], c>0$ and
   \begin{small}
     \begin{equation}\label{D4}
        \begin{aligned}
    \mathbb{D}_4=\{ x|(\frac{3}{2})^{\frac{2}{3}}(1-|\bar{\lambda}|)\zeta^{\frac{2}{9}}\in (-\infty, 0), |\bar{\lambda}|= \zeta^{-\frac{2}{3}}n^2(\pi - 2\text{arccos}\sqrt{x})^2, \zeta\to \infty, n\to \infty \},
        \end{aligned}
    \end{equation}
    \end{small}
 where  $b, \theta_3, \theta_4$ is given by (\ref{b}), (\ref{t3}), (\ref{t4}) respectively and $\bar{d}_{11}, \bar{d}_{12}, \bar{d}_{21},
     \bar{d}_{22}$ can be obtained by changing $\beta \to \alpha$ in (\ref{d11})-(\ref{d22}).
  \end{theorem}

  \subsection{Universality}

      Here we state asymptotic  of  the correlation kernel in the bulk of the spectrum and at both side of the hard-edge.

\begin{theorem}   The kernel  $K_n(x,y)$   associated   with   weight   (\ref{wpj2})  admits properties
  \begin{itemize}
      \item[(1)]When  $y \in (0,1)$, it is easy to get
      \begin{equation}
        \begin{aligned}
            \frac{1}{n}K_n(y,y) =\frac{1}{\pi \sqrt{y(1-y)}}+\mathcal{O}(\frac{1}{n}), \qquad n \to \infty,
        \end{aligned}
    \end{equation}
    and the error term is uniformly in compact subsets of $(0,1)$.
      \item[(2)] For $a \in (0,1), u,v\in \mathbb{R}$, as $n \to \infty$,
  \begin{equation}
    \begin{aligned}
        \frac{1}{n\rho(a)} K_n(a+ \frac{u}{n\rho(a)}, a+ \frac{v}{n\rho(a)})= \frac{\text{sin}(\pi(u-v))}{\pi(u-v)}+\mathcal{O}(\frac{1}{n}),
    \end{aligned}
  \end{equation}
    where $\rho(y)=(\pi \sqrt{y(1-y)})^{-1}$, and the error term is uniformly for any $a$ and $u,v$ satisfying the above conditions.
    \item[(3)] The large-$n$ behavior of the kernel near the right edge $x=1$ is
        \begin{equation}\label{kn11}
            \begin{aligned}
            \frac{1}{4n^2}K_n(1-\frac{u}{4n^2},1-\frac{v}{4n^2})=K_\Psi(u,v;2n^2t)+\mathcal{O}(\frac{1}{n^2}),
            \end{aligned}
         \end{equation}
         as $n \to \infty$ with the uniform error term for $u,v \in (0, \infty)$ and $t \in (0,c]$ where $c > 0$, in which the $\Psi$-kernel is defined as
         \begin{equation}\label{kp1}
            \begin{aligned}
                K_\Psi(u,v;\zeta) =  \frac{\psi_1(-v,\zeta)\psi_2(-u,\zeta)-\psi_1(-u,\zeta)\psi_2(-v,\zeta)}{2 \pi i(u-v)},
            \end{aligned}
         \end{equation}
        where the $\psi$-function is given by (\ref{psi121}).
        Also it is noticeable that the $\Psi$-kernel is reduced to the Bessel kernel when $\zeta \to 0^{+}$  and then the Airy kernel when $\zeta \to +\infty$. Precisely, that is
        \begin{equation}
            \begin{aligned}
                K_\Psi(u,v;\zeta) = \mathbb{J}_\beta(u,v)+\mathcal{O}(\zeta),\quad \zeta \to 0^+,
            \end{aligned}
        \end{equation}
        where the error term holds uniformly for $u,\,v$ in compact subsets of $(0,\,\infty)$ and
        \begin{equation}\label{Jbeta}
            \begin{aligned}
              \mathbb{J}_\beta(x,y)  = \frac{J_\beta(\sqrt{x})\sqrt{y}J_\beta'(\sqrt{y})-J_\beta(\sqrt{y})\sqrt{x}J_\beta'(\sqrt{x})}{2(x-y)}.
            \end{aligned}
        \end{equation}

       If the parameter $t \to 0^+$ and $n \to \infty$ so that,
       \begin{equation*}
        \begin{aligned}
            \lim_{n \to \infty}2n^2t=0,
        \end{aligned}
       \end{equation*}
       then, we can obtain the Bessel kernel limit for $K_n$
        \begin{equation}
            \begin{aligned}
                \lim_{n \to \infty}\frac{1}{4n^2}K_n(1-\frac{u}{4n^2},1-\frac{v}{4n^2})=\mathbb{J}_\beta(u,v),
            \end{aligned}
        \end{equation}
        which is uniform for $u,\,v$ in compact subsets of $(0,\infty)$.

        However, as $\zeta \to +\infty$,
        \begin{equation}\label{kpsinfity}
            \begin{aligned}
                \frac{\zeta^{4/9}}{m}K_\Psi(\zeta^{2/3}(1-\frac{u}{m\zeta^{2/9}}),\zeta^{2/3}(1-\frac{v}{m\zeta^{2/9}});\zeta) = \mathbb{A}(u,v)+\mathcal{O}(\zeta^{-2/9}),
            \end{aligned}
        \end{equation}
        in which $m=(\frac{3}{2})^{\frac{2}{3}}$ and the Airy kernel $\mathbb{A}$ is expressed by
        \begin{equation}\label{A}
            \begin{aligned}
                \mathbb{A}(x,y)=\frac{Ai(x)Ai'(y)-Ai(y)Ai'(x)}{x-y}.
            \end{aligned}
        \end{equation}
        This formula is uniform for $u,\,v$ in compact subsets of $(-\infty,\infty)$.
      Following, if the parameter $t \in (0,d]$ and $n \to \infty$ such that
        \begin{equation*}
            \begin{aligned}
                \lim_{n \to \infty}2n^2t=\infty,
            \end{aligned}
           \end{equation*}
        the Airy kernel limit for $K_n$ is given by
        \begin{equation}\label{knsinfity}
            \begin{aligned}
                \lim_{n \to \infty}\frac{s_n}{m\zeta^{2/9}} K_n(1-s_n(1-\frac{u}{m\zeta^{2/9}}),1-s_n(1-\frac{v}{m\zeta^{2/9}});\zeta)= \mathbb{A}(u,v),
            \end{aligned}
        \end{equation}
       uniformly for $u,\,v$ in compact subsets of $(-\infty,\infty)$ with $\zeta= 2n^2t,\,s_n=\zeta^{-\frac{2}{3}}/(4n^2)$.
    \item[(4)] At   left edge of the spectrum, namely near the endpoint $0$, the limit kernel has the $\Psi$-kernel asymptotic behavior as $n \to \infty$,
        \begin{equation}\label{kn01}
          \begin{aligned}
            \frac{1}{4n^2}K_n(\frac{u}{4n^2}, \frac{v}{4n^2})= \bar{K}_{\bar{\Psi}} (u,v,2n^2t) +\mathcal{O}(\frac{1}{n^2})
          \end{aligned}
        \end{equation}
        uniformly for $u,v$ in compact subsets of $(0,\infty)$ and $t$ in $(0,c]$ with a constant $c >0$. Here, the $\bar{\Psi}$-kernel is defined as
       \begin{equation}
        \begin{aligned}
            \bar{K}_{\bar{\Psi}} (u,v;\zeta)= \frac{\bar{\psi}_1(-v,\zeta) \bar{\psi}_2(-u,\zeta)-\bar{\psi}_1(-u,\zeta)\bar{\psi}_2(-v,\zeta)}{2\pi i (u-v)}
        \end{aligned}
       \end{equation}
       where the scalar function $ \bar{\psi}_(\bar{\xi},\zeta),\; k=1,2$ are given in \cite{xu2015critical} (1.22)-(1.24). In addition, from \cite{XuCritical} and \cite{xu2015critical}, by the same way as we deal with $\Psi$-kernel, we can prove the limiting $\bar{\Psi}$-kernel convert to the Bessel and Airy kernels when the parameter $\zeta$ changes from $0$ to $\infty$.
\end{itemize}
\end{theorem}

\begin{remark}
  It is also readily to show that  $K_\Psi(u,v;s)$ and $\bar{K}_{\bar{\Psi}} (u,v;s)$ are the Painlev$\acute{e}$ type kernel
   associated with the solvable model RH problem, which  admits  Lax pair,  see \cite{XuCritical} and \cite{xu2015painleve}.

\end{remark}

The  above  results can be  established  by  applying Deift-Zhou nonlinear steepest descent approach  to  the asymptotic of Pollaczek-Jacobi type OPS,
which should be first Characterized with  a appropriate  RH problem.

\section{ RH Characterization  for  the OPS}

 Following the  idea   given  by Fokas, Its and Kitaev in 1992 \cite{FokasThe},
we construct a RH problem  with respect to the weight function $w_{p_J2}(x,t)$ in (\ref{wpj2})
 as  follows.

\noindent{\bf RH problem 1.}  Find $2 \times 2$ matrix-valued function $Y(z)$ satisfying
\begin{itemize}
 \item[(a)]$Y(z)$ is analytic for $z \in  \mathbb{C} \setminus \left[ 0, 1  \right]; $
 \item[(b)]$Y(z)$ satisfies the jump condition $$Y_+(x)=Y_-(x) \begin{pmatrix}
    1 & w_{p_J2}(x,t)  \\
    0 & 1
    \end{pmatrix},	\qquad   {\rm for} \quad x \in \left( 0, 1  \right); $$
  \item[(c)]The asymptotic behavior of $Y(z)$ at infinity is
$$Y(z)=(I+\mathcal{O}(z^{-1})) \begin{pmatrix}
       z^n & 0  \\
      0 & z^{-n}
       \end{pmatrix},	\qquad  z \to  \infty; $$
  \item[(d)]The asymptotic behavior of $Y(z)$ at the origin is
$$Y(z)= \begin{pmatrix}
    \mathcal{O}(1) & \mathcal{O}(1)   \\
    \mathcal{O}(1)  & \mathcal{O}(1)
        \end{pmatrix},	\qquad  z \to  0;$$
  \item[(e)]The asymptotic behavior of $Y(z)$ at $z=1$ and $\beta >0$ is
$$Y(z)= \begin{pmatrix}
    \mathcal{O}(1) & \mathcal{O}(1)   \\
    \mathcal{O}(1)  & \mathcal{O}(1)
         \end{pmatrix},	\qquad  z \to  1. $$
\end{itemize}

In a similar way given in \cite{Deift1999Orthogonal,FokasThe},   we can show that the above RH problem for  $Y$ admits  a unique solution
  \begin{equation}
   \begin{aligned}
    Y(z)= \begin{pmatrix}
        \pi_n(z) & \frac{1}{2 \pi i } \int_{0}^{1} \frac{\pi_n(x) w_{p_J2}(x,t)}{x-z} \,\mathrm{d} x        \\
        -2 \pi i \gamma_{n-1}^2 \pi_{n-1}(z)  & -\gamma_{n-1}^2 \int_{0}^{1} \frac{\pi_{n-1}(x) w_{p_J2}(x,t)}{x-z} \,\mathrm{d} x
        \end{pmatrix},
   \end{aligned}
  \end{equation}
where $\pi_n(z) $ is the monic polynomial, and $P_{n}(z)=\gamma_{n}\pi_n(z)$ is the orthonormal polynomial with respect to the weight $ w_{p_J2}(x,t)$.

\section{Asymptotic Analysis on the RH Problem}

In this section, we will  transform    the original RH problem  for $Y$  to a solvable RH problem  via a series of invertible deformations $Y \to T \to S \to R $.
 \begin{itemize}

 \item   The first  deformation:  $Y \to T$  means that we  change the   RH problem for  $Y(z)$  to  a  normalized one, such that   $T\sim I, \ z \to \infty$.

   \item  The second  deformation:  $T \to S$
   is that we  obtain solvable RH problem for $S$   by   removing   the oscillation factors  in  RH problem for  $T$.

   \item  The third  deformation:  $S \to R$ is  to makes a small norm  RH problem for $R$  with the jumps close to the identity matrix by making of
   global parametrix $P^{(\infty)}(z)$ and local parametrixs $P^{(0)}(z)$ and $P^{(1)}(z)$ near the endpoints $z=0$ and $z=1$.
\end{itemize}

\subsection{The first deformation $Y \to T $}
In order to transform the matrix function $Y(z)$ at infinity into an identity matrix, we introduce the  following transformation
  \begin{equation}\label{T}
    \begin{aligned}
        T(z)=4^{n\sigma_3}Y(z)e^{-\frac{t}{2z(1-z)}\sigma_3}\varphi(z)^{-n\sigma_3}, \ \ z\in  \mathbb{C} \setminus \left[ 0, 1  \right],
    \end{aligned}
  \end{equation}
  where $\sigma_3$ is the Pauli matrix, $Y(z)$ is the unique solution of the above RH problem 1, and
    $\varphi(z)=2z-1+2\sqrt{z(z-1)}$ is a conformal map from $\mathbb{C} \setminus \left[ 0, 1  \right]$ onto the exterior of the unit circle.
     It can be shown that $\varphi(z)$ admits the property
     $$\varphi(z) \sim 4z, \ {\rm for} \ z\to \infty; \ \  \ \varphi_+(x)\varphi_-(x)=1, \ {\rm for } \ x \in (0,1).$$

Under the transformation (\ref{T}),  we can check  that  $T$  solves    the following RH problem.

\noindent{\bf RH problem 2.}    Find $2 \times 2$ matrix-valued function $T(z)$ satisfying
\begin{itemize}
    \item[(a)]$T(z)$ is analytic for $z \in  \mathbb{C} \setminus \left[ 0, 1  \right]. $
    \item[(b)]$T(z)$ satisfies the jump condition
    \begin{equation}\label{jump}
        \begin{aligned}
            T_+(x)=T_-(x)
           \begin{pmatrix}
             \varphi_+(x)^{-2n} & w(x)  \\
            0 & \varphi_-(x)^{-2n}
          \end{pmatrix},	\quad   {\rm for} \; x \in \left( 0, 1  \right),
        \end{aligned}
    \end{equation}
    where $ w(x)=x^{\alpha}(1-x)^{\beta}, \alpha>0, \beta>0.$
     \item[(c)]The asymptotic behavior of $T(z)$ at infinity is
     \begin{equation}
        \begin{aligned}
            T(z)=I+\mathcal{O}(z^{-1}),	\qquad  z \to  \infty
        \end{aligned}
    \end{equation}
     \item[(d)]The asymptotic behavior of $T(z)$ at the origin is
     \begin{equation}
        \begin{aligned}
            T(z)= \begin{pmatrix}
                \mathcal{O}(1) & \mathcal{O}(1)   \\
                \mathcal{O}(1)  & \mathcal{O}(1)
              \end{pmatrix}e^{-\frac{t}{2z}\sigma_3},	\qquad  z \to  0.
        \end{aligned}
    \end{equation}
     \item[(e)]The asymptotic behavior of $T(z)$ at $z=1$  is
     \begin{equation}
        \begin{aligned}
            T(z)= \begin{pmatrix}
                \mathcal{O}(1) & \mathcal{O}(1)   \\
                \mathcal{O}(1)  & \mathcal{O}(1)
             \end{pmatrix}e^{-\frac{t}{2(1-z)}\sigma_3},	\qquad  z \to  1.
        \end{aligned}
    \end{equation}
   \end{itemize}

  It is noteworthy that the transformation contains the factor $e^{-\frac{t}{2z(1-z)}\sigma_3}$ which will make $T(z)$ has two singularities at the point of $z=1$ and $z=1$ respectively. We use the similar method in \cite{xu2015painleve} to construct a new model RH problem for representing the property of singular point.

\subsection{The second transformation $T \to S $}

The jump matrix  in (\ref{jump}) can be factorized in the form
\begin{equation*}
    \begin{aligned}
        \begin{pmatrix}
         \varphi_+^{-2n} & w\\
         0 & \varphi_-^{-2n}
        \end{pmatrix}= \begin{pmatrix}
            1 & 0\\
            \varphi_-^{-2n} w^{-1 }& 1
           \end{pmatrix}  \begin{pmatrix}
            0 & w\\
            -  w^{-1 }& 0
           \end{pmatrix}\begin{pmatrix}
            1 & 0\\
            \varphi_+^{-2n} w^{-1 }& 1
           \end{pmatrix},
    \end{aligned}
\end{equation*}
which contains  two oscillation factors $\varphi_+^{-2n}$ and  $\varphi_-^{-2n}$. To remove these  oscillation factors,
 according to above factorization, we define   the  transformation
 \begin{small}
 \begin{equation}\label{S}
    \begin{aligned}
        S(z)=
           \begin{cases}
             T(z),& \text{for $z$ outside the lens region;}\\
             T(z)\begin{pmatrix}
                  1 & 0   \\
                  -w(z)^{-1}\varphi(z)^{-2n}\  & 1
                 \end{pmatrix},& \text{for $z$ in the upper lens region;}\\
             T(z)\begin{pmatrix}
                  1 & 0   \\
                  w(z)^{-1}\varphi(z)^{-2n}\  & 1
                 \end{pmatrix},& \text{for $z$ in the lower lens region,}
           \end{cases}
    \end{aligned}
 \end{equation}
\end{small}
where arg $ z\in (-\pi, \pi)$. Then $S$ solves the problem:

\noindent{\bf RH problem 3.}    Find $2 \times 2$ matrix-valued function $S(z)$ satisfying
\begin{figure}
   \begin{center}
    \begin{tikzpicture}
    \begin{scope}[line width=1pt]
    \draw[->,>=stealth] (-5,0) -- (0,0);
    \draw[-] (-0.5,0) --(5,0);
    \draw[->,>=stealth] (-3,0) arc (180:90:3 and 1.5);
    \draw (-0.1,1.5) arc (90:0:3.045 and 1.5) ;
    \draw[->,>=stealth] (-3,0) arc (-180:-90:3 and 1.5);
    \draw (-0.1,-1.5) arc (-90:0:3.045 and 1.5);
    \node[below,font=\Large] at (3.1,-0.15) {$1$};
    \node[below,font=\Large] at (-3.1,-0.15) {$0$};
    \node[above,font=\Large] at (-0.2,-1.5) {$\Sigma_3$};
    \node[above,font=\Large] at (-0.2,0.1) {$\Sigma_2$};
    \node[above,font=\Large] at (-0.2,1.5) {$\Sigma_1$};
    \end{scope}
    \end{tikzpicture}
   \end{center}
   \caption{   The contour for the RH problem for $\Sigma_S(z)$}
   \label{fig1}
\end{figure}
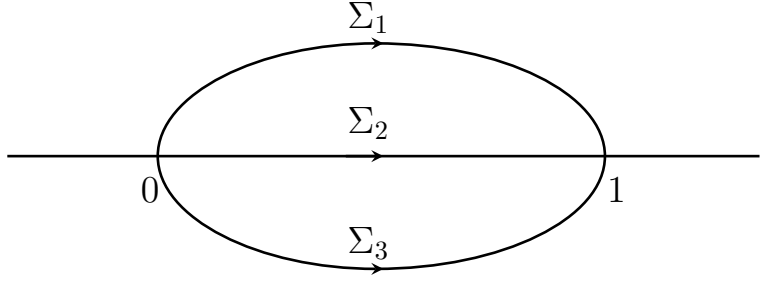

\begin{itemize}
    \item[(a)]$S(z)$ is analytic for $z \in  \mathbb{C} \setminus \bigcup_{j = 1}^{3}\Sigma_j$, illustrated in  Figure \ref{fig1}.
    \item[(b)]$S_+(z)=S_-(z)Js(z) $ for $z \in {\bigcup_{j = 1}^{3}\Sigma_j}$, and the jump $Js$ is given by
  \begin{equation}\label{sss}
    \begin{aligned}
     Js(z) &= \begin{pmatrix}
       1 & 0  \\
       w(z)^{-1}\varphi(z)^{-2n} & 1
           \end{pmatrix},   \qquad   {\rm for}\quad z \in \Sigma_1\cup \Sigma_3, \\
     Js(x) &= \begin{pmatrix}
        0& w(x)  \\
        w(x)^{-1} & 0
        \end{pmatrix},  \qquad   {\rm for}\quad x \in \Sigma_2,
    \end{aligned}
  \end{equation}
     \item[(c)]When $z \to \infty $, we get the asymptotic behavior of $S(z)$ is
   $$S(z)=I+\mathcal{O}(z^{-1}),    \qquad  z \to  \infty $$
     \item[(d)]As $z \to 0 $, the asymptotic behavior of $S(z)$ is sector wise as follows:
    \begin{small}
     \begin{equation}\label{ss0}
     S(z)=\mathcal{O}(1)e^{-\frac{t}{2z} \sigma_3}
     \left\{    \begin{aligned}
                  I,\qquad\qquad\qquad\qquad\qquad& \quad \text{outside the lens region;} \\
                  \begin{pmatrix}
                          1 & 0   \\
                         -w(z)^{-1}\varphi(z)^{-2n}\  & 1
                  \end{pmatrix},        & \quad \text{in the upper lens region;} \\
                  \begin{pmatrix}
                          1 & 0   \\
                          w(z)^{-1}\varphi(z)^{-2n}\  & 1
                          \end{pmatrix}, & \quad \text{in the lower lens region.} \\
                 \end{aligned}    \right.
     \end{equation}
    \end{small}
     \item[(e)]As $z \to 1 $, the asymptotic behavior of $S(z)$ is sector wise as follows:
     \begin{small}
      \begin{equation}\label{ss}
       S(z)=\mathcal{O}(1)e^{-\frac{t}{2(1-z)} \sigma_3}
          \left\{    \begin{aligned}
                          I,\qquad\qquad\qquad\qquad\qquad & \quad \text{outside the lens region;} \\
                         \begin{pmatrix}
                           1 & 0   \\
                           -w(z)^{-1}\varphi(z)^{-2n}\  & 1
                         \end{pmatrix}, & \quad \text{in the upper lens region;} \\
                         \begin{pmatrix}
                          1 & 0   \\
                          w(z)^{-1}\varphi(z)^{-2n}\  & 1
                         \end{pmatrix}, & \quad \text{in the lower lens region.} \\
            \end{aligned}    \right.
      \end{equation}
    \end{small}
   \end{itemize}
\subsection{Global parametrix}
When $n \to \infty$, we see that the jump matrix $Js$ on $\Sigma_1$ and $\Sigma_3$ converges to the identity matrix with an exponentially small error terms. Neglecting the exponentially small terms, $S(z)$ can be approximated by a solution of the following RH problem for $P^{(\infty)}(z)$:

\noindent{\bf RH problem 4.}    Find $2 \times 2$ matrix-valued function $P^{(\infty)}(z)$ satisfying
\begin{itemize}
    \item[(a)]$P^{(\infty)}(z)$ is analytic for $z \in  \mathbb{C} \setminus \left[ 0, 1  \right]. $
    \item[(b)]$P^{(\infty)}(z)$ satisfies the jump condition $$P^{(\infty)}_+(x)=P^{(\infty)}_-(x) \begin{pmatrix}
       0 & w(x)  \\
       -w(x)^{-1} & 0
       \end{pmatrix},	\qquad   {\rm for} \quad x \in \left( 0, 1  \right). $$
     \item[(c)]The asymptotic behavior of $P^{(\infty)}(z)$ at infinity is
     $$P^{(\infty)}(z)=(I+\mathcal{O}(z^{-1})) ,	\qquad  z \to  \infty. $$
\end{itemize}
\par We can  construct a solution to the above RH problem in the outside region applying the method in \cite{kuijlaars2002universality},
  \begin{equation}\label{pinfty}
    \begin{aligned}
        P^{(\infty)}(z)=D(\infty)^{-\sigma_3}M^{-1}a(z)^{-\sigma_3}MD(z)^{-\sigma_3},
    \end{aligned}
  \end{equation}
where $M=\frac{I+i\sigma_1}{\sqrt{2}}$, $a(z)=(\frac{z-1}{z})^{1/4}$ with arg $z \in (-\pi,\pi)$ and arg $(z-1) \in (-\pi,\pi)$, and
\begin{equation}\label{D}
    \begin{aligned}
        D(z)=\varphi^{\frac{\alpha+\beta}{2}}z^{-\frac{\alpha}{2}}(z-1)^{-\frac{\beta}{2}}
    \end{aligned}
\end{equation}
is the Szeg$\ddot{o}$ function associated with the weight $w(z)$, which is a nonzero analytic function on $\mathbb{C} \setminus \left[ 0, 1  \right]$.
Besides it is easy to verify $D_+(x)D_-(x)=w(x)$, for $x \in (0,1)$. $D(\infty)= {\lim_{z \to +\infty}}D(z) =2^{\alpha+\beta}$.
   \par After same computations,
    \begin{equation}\label{spi}
        \begin{aligned}
            S(z)P^{(\infty)}(z)^{-1}= I + \mathcal{O}(e^{-cn}),  \quad n \to \infty,
        \end{aligned}
    \end{equation}
where the error term is uniformly for $z$ away from the end-points $0$ and $1$. However, the jump matrices of $S(z)P^{(\infty)}(z)^{-1} $ are not uniformly close to the unit matrix because of  the factor $a(z)$ which has singularity. Thus, it is necessary to construct the local parametrix in the neighborhoods of the end-points.

\subsection{Local parametrix at $P^{(1)}(z)$ at $z=1$}
Firstly, we construct rhe local parametrix $P^{(1)}(z)$ in $U(1,r)=\lbrack z\in \mathbb{C} \mid z-1 \mid < r\rbrack $ for a small and fixed $r$, which solves the following RH problem:

\noindent{\bf RH problem 5.}    Find $2 \times 2$ matrix-valued function $U(1,r)$ satisfying
\begin{itemize}
    \item[(a)]$P^{(1)}(z)$ is analytic for  $z \in  U(1,r) \setminus \{\bigcup_{j = 1}^{3}\Sigma_j\}$, $\Sigma_j$ for $j=1, 2, 3$ is in Figure \ref{fig1}.
    \item[(b)]$P^{(1)}(z)$ satisfies the same jump condition with $S(z)$ on $ U(1,r)\cap\Sigma_j, j=1, 2, 3. $
    \item[(c)]$P^{(1)}(z)$ fulfills the following matching condition on $\partial U(1,r)$ and $n \to \infty$
        \begin{equation}\label{p1pi}
            \begin{aligned}
                P^{(1)}(z)P^{(\infty)}(z)^{-1}=I+\mathcal{O}(n^{-2/3}).
            \end{aligned}
        \end{equation}
    \item[(d])The asymptotic behavior of $P^{(1)}(z)$ at $z=1$ is the same as $S(z)$ in (\ref{ss}).
\end{itemize}
\par For converting all the jumps of the $P^{(1)}(z)$ to constant jumps, we apply a transformation by defining
  \begin{equation}\label{p11}
      \begin{aligned}
        P^{(1)}(z)=\hat{P}^{(1)}(z)W(z)^{-\sigma_3}\varphi(z)^{-n\sigma_3},
      \end{aligned}
  \end{equation}
where the function
   \begin{equation}\label{W}
       \begin{aligned}
          W(z)=\begin{cases}
            e^{\frac{\beta \pi i}{2}}w(z)^{\frac{1}{2}},& \text{for}\; \text{Im}z>0,\\
            e^{-\frac{\beta \pi i}{2}}w(z)^{\frac{1}{2}},& \text{for}\; \text{Im}z<0,
         \end{cases}
       \end{aligned}
   \end{equation}
 where $ W_+(x)W_-(x)=w(x)$ for $x \in (0,1)$.
Thus, it is readily seen that $\hat{P}^{(1)}(z)$ satisfies the RH problem:

\noindent{\bf RH problem 6.}    Find $2 \times 2$ matrix-valued function $\hat{P}^{(1)}(z)$ satisfying
\begin{itemize}
    \item[(a)]$\hat{P}^{(1)}(z)$ is analytic for  $z \in  U(1,r) \setminus \{\bigcup_{j = 1}^{3}\Sigma_j \}$, see Figure \ref{fig1}.
    \item[(b)]$\hat{P}^{(1)}(z)$ satisfies the jump conditions:
    \begin{equation}
        \begin{aligned}
            \hat{P}^{(1)}_+(z)=\hat{P}^{(1)}_-(z)
                                \begin{cases}
                                   \begin{pmatrix}
                                      1 & 0   \\
                                      e^{\beta \pi i}\  & 1
                                   \end{pmatrix},& \text{$z \in U(1,r)\cap \Sigma_1$,}\\
                                   \begin{pmatrix}
                                       0 & 1  \\
                                       -1\  & 0
                                  \end{pmatrix},& \text{$z \in U(1,r)\cap (0,1)$,}\\
                                  \begin{pmatrix}
                                      1 & 0   \\
                                      e^{-\beta \pi i}\  & 1
                                 \end{pmatrix},& \text{$z \in U(1,r)\cap \Sigma_3$.}
                                \end{cases}
        \end{aligned}
    \end{equation}
    \item[(c)]The asymptotic behavior at the center $z=1$ and $\beta>0$,
   \begin{equation}
       \begin{aligned}
        \hat{P}^{(1)}(z)&=\mathcal{O}(1)e^{-\frac{t}{2(1-z)}\sigma_3}W(z)^{\sigma_3}\varphi(z)^{n\sigma_3}\\
                        & \begin{cases}
                           I,& \text{outside the lens shaped region;}\\
                          \begin{pmatrix}
                            1 & 0   \\
                            -e^{\beta \pi i}\  & 1
                          \end{pmatrix},& \text{in the upper lens region;}\\
                          \begin{pmatrix}
                            1 & 0   \\
                            e^{-\beta \pi i}\  & 1
                          \end{pmatrix},& \text{in the lower lens region.}
                         \end{cases}
       \end{aligned}
   \end{equation}
\end{itemize}

 \begin{figure}
\begin{center}
    \begin{tikzpicture}
       \begin{scope}[line width=1pt]
        \draw[->,>=stealth] (-3,0) -- (-1.2,0);
        \draw[-] (-1.65,0) -- (0,0);
        \draw[dashed] (0,0) -- (3,0);
        \draw[->,>=stealth] (-2, 3) -- (-1,1.5);
        \draw[-](-1.2,1.8) -- (0,0);
        \draw[->,>=stealth] (-2, -3) -- (-1,-1.5);
        \draw[-](-1.2,-1.8) -- (0,0);
        \node[above,font=\Large] at (-1.2,0.1) {$\gamma_2$};
        \node[above,font=\Large] at (-1,1.6) {$\gamma_1$};
        \node[below,font=\Large] at (-0.8,-1.6) {$\gamma_3$};
        \node[below,font=\Large] at (0.1,-0.1) {$0$};
        \node[font=\Large] (p1) at (1.6,1.6) {$\Omega_1$};
        \node[font=\Large] (p2) at (-2,1.6) {$\Omega_2$};
        \node[font=\Large] (p3) at (-2,-1.6) {$\Omega_3$};
        \node[font=\Large] (p4) at (1.6,-1.6) {$\Omega_4$};
       \end{scope}
     \end{tikzpicture}
   \end{center}
   \caption{   Contours and regions for the model RH problem for $\Psi$.}
   \label{fig2}
\end{figure}
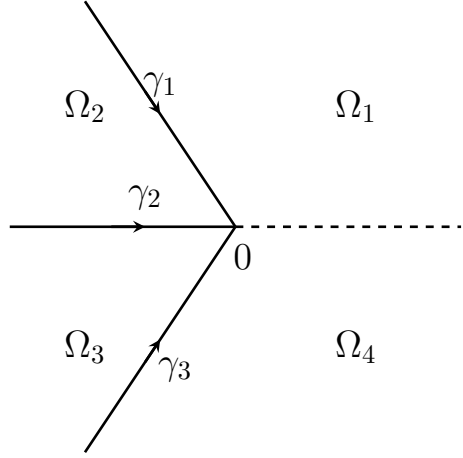

\par In order to express the above RH problem to a RH problem whose analyticity and other properties have known, using the conformal mapping, we first build a crucial model RH problem for $\Psi(\xi)=\Psi(\xi,s)$ to tackle the problem of OPS with singularities in the weight, which is similar to the model given by Xu, Dai, and Zhao \cite{XuCritical}.

\noindent{\bf RH problem 7.}    Find $2 \times 2$ matrix-valued function $\Psi(\xi)$ satisfying
\begin{itemize}
   \item[(a)]$\Psi(\xi)$ is analytic for  $\xi \in  \mathbb{C} \setminus \bigcup_{j = 1}^{3}\gamma_j. $, where $\gamma_j$ are illustrated in Figure \ref{fig2}.
   \item[(b)]$\Psi(\xi)$ satisfies the following jump conditions:

   \begin{equation}
     \begin{aligned}
        \Psi_+(\xi) &= \Psi_-(\xi)  \begin{pmatrix}
                        1 & 0  \\
                        e^{\pm \beta \pi i} & 1
        \end{pmatrix},\quad + \;\text{as} \;\xi \in  \gamma_1\; and - \text{as}  \;\xi \in  \gamma_3, \\
     \end{aligned}
   \end{equation}
   \begin{equation}
    \begin{aligned}
        \Psi_+(\xi) &= \Psi_-(\xi) \begin{pmatrix}
        0& 1  \\
        -1 & 0
        \end{pmatrix},  \qquad  \text{for}\; x \in \gamma_2,
    \end{aligned}
 \end{equation}
   \item[(c)]As $\xi \to \infty$, the asymptotic behavior of $\Psi(\xi) $ is
   \begin{equation}
    \begin{aligned}
        \Psi(\xi, s)=(I + \frac{C_1(s)}{\xi}+\mathcal{O}(\frac{1}{\xi^2}))\xi^{-\frac{1}{4}\sigma_3}\frac{I+i\sigma_1}{\sqrt{2}}e^{\sqrt{\xi}\sigma_3},
    \end{aligned}
   \end{equation}
where $$\text{arg}\; \xi \in (- \pi, \pi),\sigma_1=\begin{pmatrix}
           0 & 1 \\
           1 & 0
      \end{pmatrix},C_1(s)=\begin{pmatrix}
        q(s) & -ir(s) \\
        it(s) & -q(s)
   \end{pmatrix}$$ and from Proposition 1 in \cite{XuCritical}, $q(s)$ and $t(s)$ are expressed as a combination of $r(s)$ and $r'(s)$, respectively.
    \item[(d)] The asymptotic behavior of $\Psi(\xi)$ at $\xi =0$ is
   \begin{equation}
     \begin{aligned}
        \Psi(\xi, s)=Q(s)(I+ \mathcal{O}(\xi))e^{\frac{s}{\xi}\sigma_3}\xi^{\frac{\beta}{2}\sigma_3}
                      \begin{cases}
                        I,& \xi \in \Omega_1 \cup \Omega_4;\\
                         \begin{pmatrix}
                            1 & 0   \\
                            -e^{\beta \pi i}\  & 1
                         \end{pmatrix},& \xi \in \Omega_2 ;\\
                         \begin{pmatrix}
                             1 & 0   \\
                             e^{-\beta \pi i}\  & 1
                         \end{pmatrix},& \xi \in \Omega_3 .
                     \end{cases}
     \end{aligned}
   \end{equation}
    where arg $ \xi \in (-\pi, \pi)$ and the explicit formula for $Q(s)$ can be given by $(4,4)$ in \cite{xu2015painleve}.
 \end{itemize}
 \par The solvability of the model RH problem can be given by the similar method in \cite{XuCritical}.
     What's more, it has been proved by the Proposition 3 in \cite{XuCritical} that for $\beta>0$,
         $$r(s)=\frac{1-4\beta^2}{8}+\frac{s}{\beta}+\mathcal{O}(s^2),\qquad s \to 0,$$
         $$r(s)=\frac{3}{2}s^{2/3}-\beta s^{\frac{1}{3}}+\mathcal{O}(1),\qquad s \to +\infty.$$

    \par After comparing the RH problem for $\hat{P}^{(1)}$ with the above RH problem for $\Psi(\xi, s)$, we construct a solution $\hat{P}^{(1)}$ using the conformal mapping as follows,
     \begin{equation}
         \begin{aligned}
            \hat{P}^{(1)}(z)= E_1(z)\Psi(n^2f_1(z),2n^2t)
         \end{aligned}
     \end{equation}
    and then, by the formula (\ref{p11}) we get
    \begin{equation}\label{p1}
        \begin{aligned}
           P^{(1)}(z)= E_1(z)\Psi(n^2f_1(z),2n^2t)W(z)^{-\sigma_3}\varphi(z)^{-n\sigma_3},
        \end{aligned}
    \end{equation}
    where
    \begin{equation}\label{E1}
        \begin{aligned}
            E_1(z)=P^{(\infty)}(z)W(z)^{\sigma_3}M^{-1}(n^2f_1(z))^{\frac{1}{4}\sigma_3},
        \end{aligned}
    \end{equation}
    and $f_1(z)=(\text{log}\varphi(z))^2$.

\begin{remark}
    It is readily to verify that $f_1(z)$ is analytic in $(0,1)$, and $f_1(1)=0$. $f_1(z)$ is a conformal mapping in the neighborhood of $z=1$ and as $z \to 1$,
    $$f_1(z)=4(z-1)(1-\frac{1}{3}(z-1)+\mathcal{O}(z-1)^2).$$
    \par We transfer the RH problem for $P^{(1)}(z)$ in the $z$-plane to a model RH problem in the $\xi$-plane by the conformal mapping $\xi = n^2f_1(z) $.
    Besides, we can prove that $E_1(z)$ is analytic for $z \in U(1,r)$ owning a weak singularity at $z=1$, and by (\ref{p1}) and (\ref{pinfty}), we know
    \begin{equation}
        \begin{aligned}
            P^{(1)}(z)P^{(\infty)}(z)^{-1}=I +\mathcal{O}(n^{-1}).
        \end{aligned}
    \end{equation}
    Meanwhile, the error terms is uniformly for $z \in \partial U(1,r)\setminus \Sigma_j, j=1, 2, 3.$
\end{remark}

    \par In this situation, when $z \to 1$,
    \begin{equation}
        \begin{aligned}
            \xi = n^2 f_1(z) \approx 4 n^2 z,\; \text{and}\; \zeta = 2 n^2 t.
        \end{aligned}
    \end{equation}
    So we can apply the behavior of $\Psi(\xi,\zeta)$ as $\zeta\to 0 \; \text{and}\; \zeta \to \infty$, which have been found in \cite{ChenThe} and \cite{XuCritical}, to study the asymptotic behavior of $P^{(1)}(z)$, but we use the different parametrix $\beta$ to replace the parametrix $\alpha$ and different scaling process. For the convenience of the later part of the article, we prove the behavior briefly.
\subsubsection{The approximation of $\Psi(\xi, \zeta)$ for $\zeta \to 0 $}
\begin{figure}
\begin{center}
    \begin{tikzpicture}
    \begin{scope}[line width=1pt]
    \draw[->,>=stealth] (-3,0.09)--(-1.5,0.09);
    \draw[-] (-3,0.09)--(0,0.09);
    \draw[->,>=stealth] (0,-0.09)--(-1.5,-0.09);
    \draw[->,>=stealth] (0.1,3.15)--(-0.1,3.15);
    \draw[->,>=stealth] (-0.1,-2.85)--(0.1,-2.85);
    \draw[-] (-1.25,-0.09)--(-3,-0.09);
    \node[below,font=\Large] at (0.2,0.3) {$0$};
    \node[below,font=\Large] at (0,2.9) {$\Gamma$};
    \draw (-3,-0.1) arc (-175:181:3);
    \end{scope}
    \end{tikzpicture}
\end{center}
   \caption{    The integration path for $\Gamma$.}
   \label{fig3}
\end{figure}
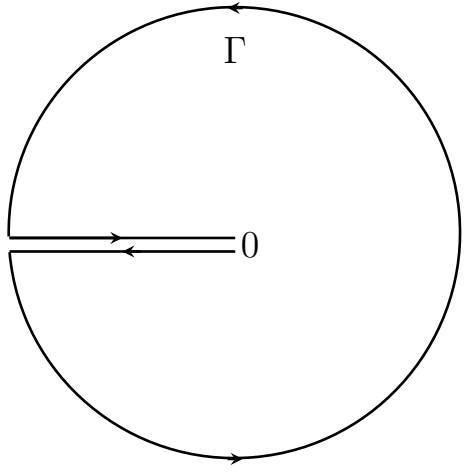

 As we can see, the essential singularity at the point $z=0$ vanishes when $\zeta=0$. Thus, the singular factor $e^{\frac{\zeta}{\xi}\sigma_3}$ can be ignored for small $\zeta$. In this case, we construct
    \begin{equation}\label{Fxz}
        \begin{aligned}
            F(\xi,\zeta)=\pi^{\frac{1}{2}\sigma_3}\begin{pmatrix}
                \xi^{-\frac{\beta}{2}}I_\beta(\sqrt{\xi}) & \frac{i}{\pi}\xi^{\frac{\beta}{2}}A_1(\sqrt{\xi})\\
                i \pi \xi^{\frac{1-\beta}{2}}I'_\beta(\sqrt{\xi}) & -\xi^{\frac{\beta+1}{2}}A_2(\sqrt{\xi })
                \end{pmatrix} \begin{pmatrix}
                    1 & f(\xi,\zeta)\\
                    0 & 1
                    \end{pmatrix} \xi^{\frac{\beta}{2}\sigma_3}e^{\frac{\zeta}{\xi}\sigma_3}J
        \end{aligned}
    \end{equation}
where \begin{equation}
       \begin{aligned}
        J=  \begin{cases}
            I,& \xi \in \Omega_1 \cup \Omega_4;\\
             \begin{pmatrix}
                1 & 0   \\
                -e^{\beta \pi i}\  & 1
             \end{pmatrix},& \xi \in \Omega_2 ;\\
             \begin{pmatrix}
                 1 & 0   \\
                 e^{-\beta \pi i}\  & 1
             \end{pmatrix},& \xi \in \Omega_3 .
            \end{cases}
      \end{aligned}
     \end{equation}
and due to the formulas of (10.25.2), (10.27.4), and (10.31.1) in \cite{Olver2010NIST}, it is readily find if $\beta > 0, \beta \notin \mathbb{N} $,
$$ A_1(z)=\frac{\pi}{2 \text{sin} (\beta \pi)}I_{-\beta}(z), \quad A_2(z)=\frac{\pi}{2 \text{sin} (\beta \pi)}I'_{-\beta}(z),$$
then if $\beta \in \mathbb{N}$,
$$A_1(z)=\frac{1}{2}(\frac{z}{2})^{-\beta}\sum_{j=0}^{\beta-1}\frac{(-1)^j (\beta -j-1)! z^{2j}}{\Gamma(j+1)4^j}+ \frac{1}{2}(\frac{-z}{2})^{\beta}\sum_{j=0}^{\infty} \frac{(\Phi(j+1)+\Phi(j+\beta+1))z^{2j}}{\Gamma(j+1)\Gamma(j+\beta+1)4^j},$$
$$A_2(z)=A'_1(z)+\frac{(-1)^{\beta+1}}{z}I_\beta(z),$$
with $\Phi(x)=\frac{\Gamma'(x)}{\Gamma(x)}$. At the same time,
 $$f(\xi,\zeta)=-\frac{1}{4 \pi \text{sin}(\beta \pi)}\int_{\Gamma} \frac{\tau^{\beta} e^{\frac{2\zeta}{\tau}}}{\tau-\xi} \,\mathrm{d} \tau = \frac{\xi^\beta}{2 \text{sin}(\beta \pi) i}(1+\mathcal{O}(\frac{\zeta}{\xi}))),\quad \beta >0, \beta \notin \mathbb{N};$$
$$f(\xi,\zeta)=\frac{(-1)^{\beta+1}}{2 \pi^2} \int_{\Gamma} \frac{\tau^{\beta} e^{\frac{2\zeta}{\tau}} \text{ln}(\frac{\sqrt{\tau}}{2})}{\tau-\xi} \,\mathrm{d} \tau=\frac{(-1)^\beta \xi^\beta \text{ln}(\frac{\sqrt{\xi}}{2})}{i \pi}(1+\mathcal{O}(\frac{\zeta}{\xi})+\mathcal{O}(\frac{\zeta^{\beta+1}}{\xi} \text{ln}\,\zeta)),\quad \beta \in \mathbb{N},$$
as $\Gamma$ has been shown in Figure \ref{fig3}.

\par After the discussion of Section 5.1 in \cite{XuCritical}, we know $\Psi(\xi,\zeta)$ is approximated by $F(\xi,\zeta)$. More precisely,
\begin{equation}\label{R_0}
    \begin{aligned}
         R_{0}(\zeta)=\begin{cases}
            \Psi(\xi,\zeta)F(\xi,0)^{-1},& \text{for}\; |\xi| >\epsilon_0 ;\\
            \Psi(\xi,\zeta)F(\xi,\zeta)^{-1},& \text{for}\; |\xi| <\epsilon_0 ;
            \end{cases}
    \end{aligned}
\end{equation}
and the jump on the circle $|\xi|=\epsilon_0$ is
\begin{equation}
    \begin{aligned}
        J_{R_{0}}(\zeta)=\begin{cases}
            I+\mathcal{O}(\zeta)+\mathcal{O}(\zeta^{\beta+1}),& \text{for} \;\beta>0, \beta \notin \mathbb{N};\\
            I+\mathcal{O}(\zeta)+\mathcal{O}(\zeta^{\beta+1}\text{ln}\,\zeta),& \text{for}\; \beta \in \mathbb{N} ;
            \end{cases}
    \end{aligned}
\end{equation}
By the norm estimation of Cauchy operator, we obtain
\begin{equation}
    \begin{aligned}
        R_{0}(\zeta)=\begin{cases}
            I+\mathcal{O}(\zeta)+\mathcal{O}(\zeta^{\beta+1}),& \text{for} \;\beta>0, \beta \notin \mathbb{N};\\
            I+\mathcal{O}(\zeta)+\mathcal{O}(\zeta^{\beta+1}\text{ln}\,\zeta),& \text{for}\; \beta \in \mathbb{N} ;
            \end{cases}
    \end{aligned}
\end{equation}
where the error terms hold uniformly for $|\xi|=\epsilon_0$ as $\zeta \to 0.$

\begin{figure}
\begin{center}
    \begin{tikzpicture}
       \begin{scope}[line width=1pt]
        \draw[->,>=stealth] (-3,0) -- (-1.2,0);
        \draw[-] (-1.65,0) -- (1,0);
        \draw[dashed] (0,0) -- (3,0);
        \draw[->,>=stealth] (-2, 3) -- (-1,1.5);
        \draw[-](-1.2,1.8) -- (0,0);
        \draw[->,>=stealth] (-2, -3) -- (-1,-1.5);
        \draw[-](-1.2,-1.8) -- (0,0);
        \draw[dashed] (1,0) -- (-1,3);
        \draw[dashed] (1,0) -- (-1,-3);
        \node[above,font=\Large] at (-1.2,0.1) {$\gamma_2$};
        \node[below,font=\Large] at (-2,2.7) {$\hat{\gamma_1}$};
        \node[above,font=\Large] at (-1.9,-2.6) {$\hat{\gamma_3}$};
        \node[font=\Large] at (-0.3,2.6) {$\gamma_1$};
        \node[font=\Large] at (-0.3,-2.6) {$\gamma_3$};
        \node[below,font=\Large] at (1.1,-0.05) {$0$};
        \node[below,font=\Large] at (0.1,-0.05) {$-I$};
        \node[font=\Large] (p1) at (1.6,1.2) {$\Pi_1$};
        \node[font=\Large] (p2) at (-2,1.2) {$\Pi_2$};
        \node[font=\Large] (p3) at (-2,-1.2) {$\Pi_3$};
        \node[font=\Large] (p4) at (1.6,-1.2) {$\Pi_4$};
        \node[font=\Large] (p5) at (-0.3,1.2) {$\Pi_5$};
        \node[font=\Large] (p6) at (-0.3,-1.2) {$\Pi_6$};
       \end{scope}
     \end{tikzpicture}
   \end{center}
   \caption{     Contours and regions for the RH problem for $H(\lambda)$.}
   \label{fig4}
\end{figure}
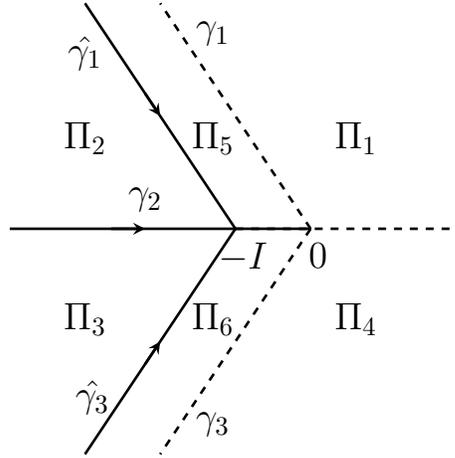

\subsubsection{The approximation of $\Psi(\xi, \zeta)$ for $\zeta \to \infty $}
In order to removal of singularity of $\Psi(\xi, \zeta)$ as $\zeta \to \infty$, we rescale $\xi$ as $\zeta^{\frac{2}{3}}\lambda$, then we have
 \begin{equation}\label{H}
    \begin{aligned}
         H(\lambda)=\zeta^{\frac{1}{6}\sigma_3}\Psi(\zeta^{\frac{2}{3}}\lambda,\zeta)e^{-\zeta^{\frac{1}{3}}\theta(\lambda)\sigma_3},\quad \theta(\lambda)=\frac{(\lambda+1)^{\frac{2}{3}}}{\lambda}.
    \end{aligned}
 \end{equation}
 which solves the following RH problem:

  \noindent{\bf RH problem 8.}    Find $2 \times 2$ matrix-valued function $H(\lambda)$ satisfying
 \begin{itemize}
    \item[(a)]$H(\lambda)$ is analytic for $\lambda \in \mathbb{C} \setminus {\bigcup_{j = 1}^{3}\gamma_j} $, see Figure \ref{fig4}.
    \item[(b)]$H(\lambda)$ fulfills the jump conditions:
    \begin{equation}
        \begin{aligned}
             H_+(\lambda)=H_-(\lambda) \begin{cases}
                \begin{pmatrix}
                    1 & 0   \\
                    e^{\beta \pi i} e^{-2\zeta^{\frac{1}{3}}\theta(\lambda)} & 1
                \end{pmatrix}, & \; \lambda \in \gamma_1,\\
             \begin{pmatrix}
                 0 & 1   \\
                 -1 & 0
             \end{pmatrix}, & \; \lambda \in (-\infty,-1),\\
             \begin{pmatrix}
                0 &   e^{2\zeta^{\frac{1}{3}}\theta(\lambda)}  \\
                - e^{-2\zeta^{\frac{1}{3}}\theta(\lambda)} & 0
            \end{pmatrix}, & \; \lambda \in (-1,0),\\
            \begin{pmatrix}
                1 & 0   \\
                e^{-\beta \pi i} e^{-2\zeta^{\frac{1}{3}}\theta(\lambda)} & 1
            \end{pmatrix}, & \; \lambda \in \gamma_3;\\
             \end{cases}
        \end{aligned}
    \end{equation}
    \item[(c)] The asymptotic behavior of $H(\lambda)$ when $\lambda \to \infty$ is
    \begin{equation}
        \begin{aligned}
          H(\lambda) = (I + \mathcal{O}(\frac{1}{\lambda}))\lambda^{-\frac{1}{4}\sigma_3}M;
        \end{aligned}
    \end{equation}
    \item[(d)] The asymptotic behavior of $H(\lambda)$ when $\lambda \to 0$ is
    \begin{equation}
        \begin{aligned}
            H(\lambda) =\mathcal{O}(1) \lambda^{\frac{\beta}{2}\sigma_3}\begin{cases}
                I,& \lambda \in \Pi_1 \cup \Pi_4;\\
                \begin{pmatrix}
                   1 & 0   \\
                   -e^{\beta \pi i}e^{-2\zeta^{\frac{1}{3}}\theta(\lambda)}\  & 1
                \end{pmatrix},& \lambda \in \Pi_2 ;\\
                \begin{pmatrix}
                    1 & 0   \\
                    e^{-\beta \pi i}e^{-2\zeta^{\frac{1}{3}}\theta(\lambda)}\  & 1
                \end{pmatrix},& \lambda \in \Pi_3 .
            \end{cases}
        \end{aligned}
    \end{equation}
 \end{itemize}
Then, to move the jumps from $\gamma_1$ and $\gamma_3$ to $\hat{\gamma_1}$ and $\hat{\gamma_3}$, we introduce the following transformation
 \begin{equation}\label{N}
    \begin{aligned}
        N(\lambda)=\begin{cases}
                 H(\lambda) \begin{pmatrix}
                              1 & 0   \\
                              e^{\beta \pi i} e^{-2\zeta^{\frac{1}{3}}\theta(\lambda)} &1
                            \end{pmatrix}, & \quad \lambda \in \Pi_5,\\
                 H(\lambda) \begin{pmatrix}
                                1 & 0   \\
                                -e^{-\beta \pi i} e^{-2\zeta^{\frac{1}{3}}\theta(\lambda)} &1
                              \end{pmatrix}, & \quad \lambda \in \Pi_6,\\
                 H(\lambda),   & \quad \text{otherwise}.
        \end{cases}
    \end{aligned}
 \end{equation}
 So, the RH problem for $ N(\lambda)$ is constructed as follows:

  \noindent{\bf RH problem 9.}    Find $2 \times 2$ matrix-valued function $N(\lambda)$ satisfying
 \begin{itemize}
    \item[(a)]$N(\lambda)$ is analytic for $\lambda \in \mathbb{C} \setminus \hat{\gamma_1}\cup \hat{\gamma_3}\cup (-\infty,-1)\cup (-1,0) $, see Figure \ref{fig4}.
    \item[(b)]$N(\lambda)$ fulfills the jump conditions:
    \begin{equation}
        \begin{aligned}
             N_+(\lambda)=N_-(\lambda) \begin{cases}
                \begin{pmatrix}
                    1 & 0   \\
                    e^{\beta \pi i} e^{-2\zeta^{\frac{1}{3}}\theta(\lambda)} & 1
                \end{pmatrix}, & \; \lambda \in \hat{\gamma_1},\\
             \begin{pmatrix}
                 0 & 1   \\
                 -1 & 0
             \end{pmatrix}, & \; \lambda \in (-\infty,-1),\\
             \begin{pmatrix}
                0 &   e^{2\zeta^{\frac{1}{3}}\theta(\lambda)}  \\
                - e^{-2\zeta^{\frac{1}{3}}\theta(\lambda)} & 0
            \end{pmatrix}, & \; \lambda \in (-1,0),\\
            \begin{pmatrix}
                1 & 0   \\
                e^{-\beta \pi i} e^{-2\zeta^{\frac{1}{3}}\theta(\lambda)} & 1
            \end{pmatrix}, & \; \lambda \in \hat{\gamma_3};\\
             \end{cases}
        \end{aligned}
    \end{equation}
    \item[(c)] The asymptotic behavior of $N(\lambda)$ when $\lambda \to \infty$ is
    \begin{equation}
        \begin{aligned}
          N(\lambda) = (I + \mathcal{O}(\frac{1}{\lambda}))\lambda^{-\frac{1}{4}\sigma_3}M;
        \end{aligned}
    \end{equation}
    \item[(d)] The asymptotic behavior of $N(\lambda)$ when $\lambda \to 0$ is
    \begin{equation}
        \begin{aligned}
            N(\lambda) =\mathcal{O}(1) \lambda^{\frac{\beta}{2}\sigma_3}.
        \end{aligned}
    \end{equation}
 \end{itemize}
\par Next, the global parametrix $U(\lambda)$ can be constructed by ignoring the exponentially small entries in the above jumps, which satisfies the following RH problem:

\noindent{\bf RH problem 10.}    Find $2 \times 2$ matrix-valued function $U(\lambda)$ satisfying
\begin{itemize}
    \item[(a)]$U(\lambda)$ is analytic for $\lambda \in \mathbb{C} \setminus (-\infty,0).$
    \item[(b)]$U(\lambda)$ fulfills the jump conditions:
    \begin{equation}
        \begin{aligned}
             U_+(\lambda)=U_-(\lambda) \begin{cases}
                e^{\beta \pi i\sigma_3}, & \; \lambda \in (-1,0),\\
             \begin{pmatrix}
                 0 & 1   \\
                 -1 & 0
             \end{pmatrix}, & \; \lambda \in (-\infty,-1);
             \end{cases}
        \end{aligned}
    \end{equation}
    \item[(c)] The asymptotic behavior of $U(\lambda)$ at infinity is the same as $N(\lambda).$
    \item[(d)] The asymptotic behavior of $U(\lambda)$ at the origin is also the same as $N(\lambda).$
 \end{itemize}
By the similar technique in \cite{its2009asymptotics}, we find
 \begin{equation}\label{U}
    \begin{aligned}
        U(\lambda)=(\lambda+1)^{-\frac{1}{4}\sigma_3}M(\frac{\sqrt{\lambda+1}+1}{\sqrt{\lambda+1}-1})^{-\frac{\beta}{2}\sigma_3},
    \end{aligned}
\end{equation}
with arg$\,(\lambda+1) \in (-\pi, \pi)$ and  arg$\,\lambda \in (-\pi, \pi).$
It is obvious that $N(\lambda)$ can be approximated by $U(\lambda)$ except for the neighborhood of the singularity at $\lambda=-1$, namely $U(-1,r)$, so we need to construct an extra model for the local parametrix, which is a Airy model RH problem. That means that the local parametrix can be expressed as
    \begin{equation}\label{U1}
        \begin{aligned}
        U_1(\lambda)=E_{01}(\lambda)\Psi_A(\zeta^{\frac{2}{9}}f_{01}(\lambda))e^{-\zeta^{\frac{1}{3}}\theta(\lambda)\sigma_3}e^{\pm \frac{1}{2}\beta \pi i}, \quad \text{for} \; \pm\text{Im}\lambda >0, \lambda \in U(-1,r),
        \end{aligned}
    \end{equation}
where $\Psi_A$  is the Airy model RH problem in section 7 of \cite{deift1999strong} and
    \begin{equation}\label{E01}
      \begin{aligned}
        E_{01}(\lambda) =U(\lambda)e^{\mp \frac{\beta \pi i}{2}\sigma_3}M^{-1}(\zeta^{\frac{2}{9}}f_{01}(\lambda))^{\frac{1}{4}\sigma_3}, \quad \text{for} \; \pm\text{Im}\lambda >0, \lambda \in U(-1,r),
      \end{aligned}
    \end{equation}
    \begin{equation}
        \begin{aligned}
            f_{01}(\lambda)=(-\frac{3}{2}\theta(\lambda))^{\frac{2}{3}}.
        \end{aligned}
    \end{equation}
\par Then, we introduce the final transformation
    \begin{equation}\label{R01}
        \begin{aligned}
            R_{01}=\begin{cases}
                N(\lambda)U(\lambda)^{-1}, & |\lambda+1|>r,\\
                N(\lambda)U_1(\lambda)^{-1}, & |\lambda+1|<r.
            \end{cases}
        \end{aligned}
    \end{equation}
and it is also easy to check that, for $\zeta \to \infty$
   \begin{equation}
     \begin{aligned}
        J_{R_{01}}=\begin{cases}
                   I +\mathcal{O}(e^{-c\zeta^{\frac{1}{3}}}), & \Sigma_{R_{01}} \setminus \partial U(-1,r),\\
                   I +\mathcal{O}(\zeta^{\frac{1}{3}}), & |\lambda+1|=r.
            \end{cases}
     \end{aligned}
   \end{equation}

    Applying the norm estimation of Cauchy operator, it follows that
\begin{equation}\label{R01p}
    \begin{aligned}
        R_{01}(\lambda)=I+ \mathcal{O}(\zeta^{-\frac{1}{3}}), \quad \zeta\to \infty.
    \end{aligned}
\end{equation}
Then, using the technique showed in section 8 of \cite{kuijlaars2001riemann} or the result in section 2.6.2 of \cite{ChenThe}, we have
\begin{equation}\label{R01pp}
    \begin{aligned}
        R_{01}(\lambda)=I+ (\frac{M_1}{\lambda+1}+\frac{M_2}{(\lambda+1)^2})\zeta^{-\frac{1}{3}}+\mathcal{O}(\zeta^{-\frac{1}{3}}), \quad \zeta\to \infty.
    \end{aligned}
\end{equation}
where $M_1=\frac{7}{72}\begin{pmatrix}
    0 & 0 \\
    i & 0
\end{pmatrix}$ and $M_2=\frac{5}{72}\begin{pmatrix}
    0 & i \\
    0 & 0
\end{pmatrix}$.

\subsection{Local parametrix at $P^{(0)}(z)$ at $z=0$}
Next,  the local parametrix $P^{(0)}(z)$ in the neighborhood of $U(0,r)=\lbrack {z\in \mathbb{C}:|z|<r} \rbrack$ with a small and fixed $r$ is constructed, which satisfies the following RH problem:

\noindent{\bf RH problem 11.}    Find $2 \times 2$ matrix-valued function $P^{(0)}(z)$ satisfying
 \begin{itemize}
   \item[(a)]$P^{(0)}(z)$ is analytic for $U(0,r) \setminus {\bigcup_{j = 1}^{3}\Sigma_j}$, where $\Sigma_j$ is in Figure \ref{fig1}.
   \item[(b)]The jump condition of $P^{(0)}(z)$ is same as $S(z)$ on $U(0,r)\cap \Sigma_j, j=1, 2, 3.$
   \item[(c)]The matching condition on the boundary $\partial U(0,r)$ is when $n \to \infty$
     \begin{equation}\label{p0pi}
         \begin{aligned}
            P^{(0)}(z)P^{(\infty)}(z)^{-1}= I +\mathcal{O}(n^{-\frac{2}{3}}),
         \end{aligned}
     \end{equation}
        and the error terms is uniformly for $z \in \partial U(0,r)\setminus \Sigma_j, j=1, 2, 3.$
    \item[(d)]The asymptotic behavior of $ P^{(0)}(z)$ at $z=0$ is also the same as $S(z)$ in (\ref{ss0}).
    \end{itemize}
    \par To transform the jumps of $P^{(0)}(z)$ into constant jumps, we introduce the transformation as follow:
     \begin{equation}\label{p0}
        \begin{aligned}
            P^{(0)}(z)=\hat{P}^{(0)}(z)V(z)^{-\sigma_3}\varphi(z)^{-n\sigma_3}
        \end{aligned}
    \end{equation}
    in which the function $V(z)$ is described by
    \begin{equation}\label{V}
        \begin{aligned}
            V(z)=\begin{cases}
                e^{-\frac{\alpha \pi i}{2}}w(z)^{\frac{1}{2}},& \text{for}\; \text{Im}z>0,\\
                e^{\frac{\alpha \pi i}{2}}w(z)^{\frac{1}{2}},& \text{for}\; \text{Im}z<0,
             \end{cases}
        \end{aligned}
    \end{equation}
    and $V_+(x)V_-(x)=w(x)$ for $x\in (0,1)$. Also, when $z \to 0$, $V(z)$ can be approximated by $(-z)^{\frac{\alpha}{2}}$.
    \par Thus, we can directly verify that $\hat{P}^{(0)}(z)$ subjects to the RH problem as follow:

    \noindent{\bf RH problem 12.}    Find $2 \times 2$ matrix-valued function $\hat{P}^{(0)}(z)$ satisfying
    \begin{itemize}
        \item[(a)]$\hat{P}^{(0)}(z)$ is analytic for $z \in  \mathbb{C} \setminus {\bigcup_{j = 1}^{3}\Sigma_j}$, see Figure \ref{fig1}.
        \item[(b)]$\hat{P}^{(0)}(z)$ satisfies the jump conditions:
        \begin{equation}
            \begin{aligned}
                \hat{P}^{(0)}_+(z)=\hat{P}^{(0)}_-(z)
                                    \begin{cases}
                                       \begin{pmatrix}
                                          1 & 0   \\
                                          e^{- \alpha \pi i}\  & 1
                                       \end{pmatrix},& \text{$z \in U(0,r)\cap \Sigma_1$,}\\
                                       \begin{pmatrix}
                                           0 & 1  \\
                                           -1\  & 0
                                      \end{pmatrix},& \text{$z \in U(0,r)\cap (0,1)$,}\\
                                      \begin{pmatrix}
                                          1 & 0   \\
                                          e^{\alpha \pi i}\  & 1
                                     \end{pmatrix},& \text{$z \in U(1,r)\cap \Sigma_3$.}
                                    \end{cases}
            \end{aligned}
        \end{equation}
        \item[(c)]The asymptotic behavior at the center $z=0$,
       \begin{equation}
          \begin{aligned}
             \hat{P}^{(0)}(z)&=\mathcal{O}(1)e^{-\frac{t}{2z}\sigma_3}V(z)^{\sigma_3}\varphi(z)^{n\sigma_3}\\
                        & \begin{cases}
                           I,& \text{outside the lens region;}\\
                          \begin{pmatrix}
                            1 & 0   \\
                            -e^{- \alpha \pi i}\  & 1
                          \end{pmatrix},& \text{in the upper lens region;}\\
                          \begin{pmatrix}
                            1 & 0   \\
                            e^{\alpha \pi i}\  & 1
                          \end{pmatrix},& \text{in the lower lens region.}
                         \end{cases}
            \end{aligned}
         \end{equation}
    \end{itemize}
    \par After given the RH problem of $\hat{P}^{(0)}(z)$, we think about a model RH problem which is introduced by Xu, Dai and Zhao \cite{XuCritical} and then developed by Chen and other experts \cite{ChenThe} to deal with the local parametrix of $\hat{P}^{(0)}(z)$.
    To prevent confusion, we record this model as $\bar{\Psi}(\bar{\xi})=\bar{\Psi}(\bar{\xi},\bar{s})$.
    \par In order to find out the conformal mapping, we first define auxiliary function $h(z)$:
    \begin{equation}
        \begin{aligned}
             h(z)=\text{log}\; \varphi(z)=\text{log}\; (2z-1+2\sqrt{z(z-1)})
        \end{aligned}
    \end{equation}
    and $h(z)$ is analytic in $\mathbb{C}\setminus (-\infty,1]$ apparently. As $\varphi_+(x)\varphi_-(x)=1$, we get $h_+(x)=-h_-(x)$. Then we define
    \begin{equation}
        \begin{aligned}
            f_0(z)=(h(z)-h(0))^2=(\text{log}\;\frac{\varphi(z)}{\varphi(0)})^2
        \end{aligned}
    \end{equation}
    and when $z\to 0$, the asymptotic behavior of $f_0(z)$ is $f_0(z)=-4z(1+\frac{1}{3}z+\mathcal{O}(z^2))$ and $f_0'(0)=-4$. Therefore $f_0(z)$ is a conformal mapping in $U(0,r)$ where $r>0$ is small and fixed.
    \par Compared the model RH problem of $\bar{\Psi}(\bar{\xi}, \bar{s})$ with the local parametrix $\hat{P}^{(0)}(z)$, the consequence is that $\hat{P}^{(0)}(z)$ can be expressed by
    \begin{equation}
        \begin{aligned}
            \hat{P}^{(0)}(z)= E_0(z)\bar{\Psi}(n^2f_0(z),2n^2t)e^{-\frac{\pi i}{2}\sigma_3}
        \end{aligned}
    \end{equation}
    This means $\bar{\xi}=n^2f_0(z)$ and $\bar{\zeta}= 2n^2 t$ which equals to $\zeta$. Therefore, in order to avoid repetition, we will use the same symbol $\zeta$ in a later article. In the meantime, by the formula (\ref{p0}),
    \begin{equation}\label{p00}
        \begin{aligned}
            P^{(0)}(z)= E_0(z)\bar{\Psi}(n^2f_0(z),2n^2t)e^{-\frac{\pi i}{2}\sigma_3}V(z)^{-\sigma_3}\varphi(z)^{-n\sigma_3}
        \end{aligned}
    \end{equation}
    where
    \begin{equation}\label{E0}
        \begin{aligned}
            E_0(z)= P^{(\infty)}(z) \varphi(0)^{n \sigma_3} V(z)^{\sigma_3} e^{\frac{\pi i}{2}\sigma_3} (\frac{I + i \sigma_1}{\sqrt{2}})^{-1} (n^2f_0(z))^{\frac{1}{4}\sigma_3}
         \end{aligned}
     \end{equation}
    \par Just as we thought about $E_1(z)$, $E_0(z)$ is analytic for $z \in  U(0,r)$ with a weak singularity at the point $z=0$, and
    \begin{equation}
        \begin{aligned}
            P^{(0)}(z)P^{(\infty)}(z)^{-1}= I +\mathcal{O}(\frac{1}{n})
        \end{aligned}
    \end{equation}
     can be achieved by the formula (\ref{pinfty}) and (\ref{p00}), where the error terms hold uniformly for $z \in \partial U(0,r)\setminus \Sigma_j, j=1, 2, 3.$
     \par  In a similar  way  to  the asymptotic properties of $P^{(1)}(z)$,   we also can study the behavior of $P^{(0)}(z)$ by using the asymptotic behavior
      of $\bar{\Psi}(\bar{\xi}, \bar{s})$ which have been researched by Chen and other scholars in \cite{ChenThe}. Here we refer directly to the results in section 2.6 in the reference \cite{ChenThe}.

\subsection{The final transformation $S \to R $}
   In the final transformation, we define $R(z)$ as follow, with the results of the global parametrix $P^{(\infty)}(z)$ and local parametrixs $P^{(0)}(z)$ and $P^{(1)}(z)$,

   \begin{equation}\label{R}
    \begin{aligned}
          R(z)& =\begin{cases}
               S(z)P^{(\infty)}(z)^{-1},& z \in \mathbb{C} \setminus U(1,r)\cap U(0,r)\cap \Sigma_S;\\
               S(z)P^{(1)}(z)^{-1},& z \in U(1,r) \setminus \Sigma_S; \\
               S(z)P^{(0)}(z)^{-1},& z \in U(0,r) \setminus \Sigma_S.
          \end{cases}
    \end{aligned}
   \end{equation}
    \par Thus, $R(z)$ is said to satisfy the following RH problem:

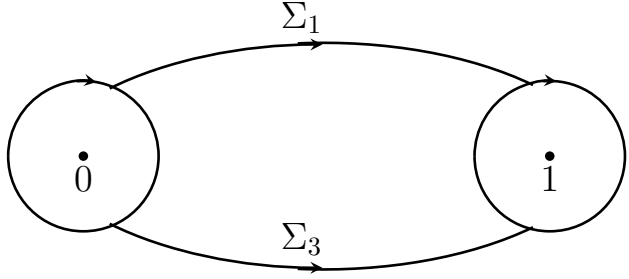
\begin{figure}
    \begin{center}
        \begin{tikzpicture}
           \begin{scope}[line width=1pt]
            \draw[->,>=stealth] (-2.9,0) arc (180:80:1);
            \draw (-2.9,0) arc (-180:95:1);
            \draw[->,>=stealth] (3.3,0) arc (180:85:1);
            \draw (3.3,0) arc (-180:95:1);
            \draw [->,>=stealth] (-1.55,0.9) arc (135:90:4 and 2);
            \draw (0.95,1.5) arc (95:46:4 and 2);
            \draw [->,>=stealth] (-1.55,-0.9) arc (-135:-90:4 and 2);
            \draw (0.95,-1.5) arc (-95:-46:4 and 2);
           \node[below,font=\Large] at (-1.9,0) {$0$};
           \fill(-1.9,0) circle(1.8pt);
           \node[below,font=\Large] at (4.3,0) {$1$};
           \fill(4.3,0) circle(1.8pt);
           \node[above,font=\Large] at (1,1.5) {$\Sigma_1$};
           \node[above,font=\Large] at (1,-1.45) {$\Sigma_3$};
           \end{scope}
        \end{tikzpicture}
      \end{center}
   \caption{     The contour for the RH problem for $R(z)$.}
   \label{fig5}
\end{figure}

    \noindent{\bf RH problem 13.}    Find $2 \times 2$ matrix-valued function $R(z)$ satisfying
    \begin{itemize}
        \item[(a)]$R(z)$ is analytic for $z \in  \mathbb{C} \setminus \Sigma_R$, as in Figure \ref{fig5}.
        \item[(b)] The jump condition of $R(z)$ is
          \begin{equation}
            \begin{aligned}
                R_+(z)=R_-(z)J_R(z),\quad z \in \Sigma_R,
            \end{aligned}
          \end{equation}
        where
           \begin{equation}\label{p0p1pi}
            \begin{aligned}
                  J_R(z) & = \begin{cases}
                   P^{(\infty)}(z)J_S(z)P^{(\infty)}(z)^{-1},& z \in \Sigma_1 \cup \Sigma_3,\\
                   P^{(0)}(z)P^{(\infty)}(z)^{-1},& z \in \partial U(0,r),\\
                   P^{(1)}(z)P^{(\infty)}(z)^{-1},& z \in \partial U(1,r).
               \end{cases}
            \end{aligned}
           \end{equation}
        and $J_S(z)$ represents the jump matrices of $S(z)$ in (\ref{sss}).
        \item[(c)] The asymptotic behavior of $R(z)$ at infty:
                 $$R(z)=I+\mathcal{O}(z^{-1}),    \qquad  z \to  \infty$$
    \end{itemize}
                 \par At the same time,  using the relations in (\ref{spi}),\; (\ref{p1pi}),\; (\ref{p0pi}) and (\ref{p0p1pi}), we find,
                 \begin{equation}
                    \begin{aligned}
                        J_R(z)= & = \begin{cases}
                            I+\mathcal{O}(e^{-cn}),& z \in \Sigma_1 \cup \Sigma_3,\\
                            I+\mathcal{O}(n^{-\frac{2}{3}}),& z \in \partial U(1,r) \cup \partial U(0,r).
                        \end{cases}
                    \end{aligned}
                 \end{equation}
                where $c$ is a positive constant, and the error terms is uniformly for $z$ on the corresponding contours. Hence, we have
                \begin{equation}\label{JR}
                    \begin{aligned}
                        ||J_R(z)-I||_{L^2\cup L^{\infty}(\Sigma_R)}= \mathcal{O}(n^{-\frac{2}{3}})
                    \end{aligned}
                \end{equation}
            Hence, applying the Plemelj formula in \cite{Deift1999Orthogonal} and \cite{deift1999strong},
            \begin{equation}\label{Rp}
                \begin{aligned}
                    R(z)=I+\frac{1}{2 \pi i} \int_{\Sigma_R} \frac{R_-(\tau) (J_R(\tau)-I)}{\tau-z} \,\mathrm{d} \tau, \quad z \notin  \Sigma_R,
                \end{aligned}
            \end{equation}
            Finally, combined (\ref{JR}) with (\ref{Rp}),we have
            \begin{equation}
                \begin{aligned}
                    R(z)=I+\mathcal{O}(n^{-\frac{2}{3}}),
                \end{aligned}
            \end{equation}
           uniformly for $z$ in the whole complex plane.
          \par So far, we have completed the Deift-Zhou steepest descent analysis.

\section{The Proofs of Main Results}
   We first start with the proof of Theorem 1 about the asymptotics for monic OPS with the two singularities Pollaczek-Jacobi type weight.
   \\ \textbf{Proof of Theorem 1.} For $z \in \mathbb{C} \setminus \left[ 0, 1  \right]$, using the formula (\ref{T}),\;(\ref{S}),\;(\ref{pinfty}),\;and (\ref{R}),
     \begin{equation*}
        \begin{aligned}
        Y(z) & =4^{-n\sigma_3}T(z)e^{\frac{t}{2z(1-z)}\sigma_3}\varphi(z)^{n\sigma_3}\\
             & =4^{-n\sigma_3}S(z)e^{\frac{t}{2z(1-z)}\sigma_3}\varphi(z)^{n\sigma_3}\\
             & =4^{-n\sigma_3}R(z)P^{(\infty)}(z)e^{\frac{t}{2z(1-z)}\sigma_3}\varphi(z)^{n\sigma_3}\\
             & =4^{-n\sigma_3}R(z)D(\infty)^{-\sigma_3}M^{-1}a(z)^{-\sigma_3}MD(z)^{-\sigma_3}(z)e^{\frac{t}{2z(1-z)}\sigma_3}\varphi(z)^{n\sigma_3} \\
         \end{aligned}
     \end{equation*}
     due to
     \begin{equation*}
        \begin{aligned}
            \begin{pmatrix}
                Y_{11}(z)\\
                Y_{22}(z)
            \end{pmatrix}= Y(z)\begin{pmatrix}
                                   1\\
                                   0
                                \end{pmatrix}
        \end{aligned}
    \end{equation*}
    Then, we get
    \begin{equation*}
        \begin{aligned}
            \begin{pmatrix}
                Y_{11}(z)\\
                Y_{22}(z)
            \end{pmatrix}= 4^{-n\sigma_3}R(z)\begin{pmatrix}
                2^{-(\alpha+\beta)} \frac{a(z)^{-1}+a(z)}{2} D(z)e^{\frac{t}{2z(1-z)}}\varphi(z)^{n}\\
                i 2^{\alpha+\beta} \frac{a(z)-a(z)^{-1}}{2} D(z)e^{\frac{t}{2z(1-z)}}\varphi(z)^{n}
                                             \end{pmatrix}.
        \end{aligned}
    \end{equation*}
    Plus, $Y_{11}(z)=\pi_n(z)$, $a(z)+a(z)^{-1}=\varphi(z)^{\frac{1}{2}}(z-1)^{-\frac{1}{4}}z^{-\frac{1}{4}}$, $D(z)=\varphi(z)^{\frac{\alpha+\beta}{2}}(z-1)^{-\frac{\beta}{2}}z^{-\frac{\alpha}{2}}$. Finally, after some calculation, the result of Theorem 1 has been obtained.
    In order to reduce complex calculations, our next proofs start from the first column of matrix $Y(z)$, which is enough for us to calculate $\pi_n(z)$.
    \par Secondly, the asymptotic behavior of $\pi_n(z)$ in the domain $(0,1)$ has been proved as follow.
\\ \textbf{Proof of Theorem 2.} $\pi_n(z)$ is a polynomial function, which is analytical and continuous in the complex plane. So it is certainly continuous on the whole axis, especially for $z \in (0,1)$. According to the definition of the continuity of complex functions,
    \begin{equation*}
     \begin{aligned}
        \pi_n(z)= \lim\limits_{\substack{z' \to z \\ z'\in \mathbb{C}}} \pi_n(z')
     \end{aligned}
    \end{equation*}
   where $z'\in \mathbb{C}$, and the limit that tends to $z$ is from any orientation in the field of $z$ in the complex plane.
   Therefore, we consider the condition that $z'$ tends to $z$ on upper plane.
   Besides, by tracing back the steps $Y \to T\to S\to R$, that is (\ref{T}),\;(\ref{S}),\;(\ref{pinfty}),\;and (\ref{R}), for $z \in (0,1)$, we obtain that
    \begin{equation*}
     \begin{aligned}
            \begin{pmatrix}
                Y_{11}(z)\\
                Y_{22}(z)
            \end{pmatrix}& =4^{-n\sigma_3}R(z)2^{-(\alpha+\beta)\sigma_3}\begin{pmatrix}
                \frac{a_+(z)^{-1}+a_+(z)}{2} & i \frac{a_+(z)^{-1}-a_+(z)}{2}\\
                i \frac{a_+(z)-a_+(z)^{-1}}{2} &  \frac{a_+(z)^{-1}+a_+(z)}{2}
            \end{pmatrix}\\
            & \begin{pmatrix}
                D_+(z) \varphi_+(z)^n e^{\frac{t}{2z(1-z)}} \\
               D_+(z)^{-1} z^{-\alpha} (1-z)^{-\beta} \varphi_+(z)^{-n} e^{\frac{t}{2z(1-z)}}
             \end{pmatrix}
     \end{aligned}
    \end{equation*}
    By the formulas (\ref{pinfty}) and (\ref{D}), it is easily to check that
    \begin{equation}
        \begin{aligned}
            D_+(z)=e^{i((\alpha+\beta)\text{arccos}\sqrt{z}-\beta \pi / 2)} z^{-\frac{\alpha}{2}}(1-z)^{-\frac{\beta}{2}}
        \end{aligned}
    \end{equation}
    \begin{equation}
        \begin{aligned}
        \varphi_+(z)^n=e^{i 2 n \text{arccos}\sqrt{z}}
        \end{aligned}
    \end{equation}
    \begin{equation}
        \begin{aligned}
            a_+(z)^{-1}+a_+(z)& =\varphi_+(z)^{\frac{1}{2}}(z-1)_+^{-\frac{1}{4}}z_+^{-\frac{1}{4}}\\
                              & =e^{i(\text{arccos}\sqrt{z}-\frac{\pi}{4})}(1-z)^{-\frac{1}{4}}z^{-\frac{1}{4}}
        \end{aligned}
    \end{equation}
    \begin{equation}
        \begin{aligned}
            i(a_+(z)^{-1}-a_+(z))=e^{i(-\text{arccos}\sqrt{z}+\frac{\pi}{4})}(1-z)^{-\frac{1}{4}}z^{-\frac{1}{4}}
        \end{aligned}
    \end{equation}
    then the above matrix turns into
    \begin{equation*}
        \begin{aligned}
               \begin{pmatrix}
                   Y_{11}(z)\\
                   Y_{22}(z)
               \end{pmatrix} & = (w_{p_J2}(z,t))^{-\frac{1}{2}} z^{-\frac{1}{4}} (1-z)^{-\frac{1}{4}} 4^{-n\sigma_3} R(z) 2^{-(\alpha+\beta)\sigma_3}\\
                             & \begin{pmatrix}
                               \text{cos}((2n+\alpha+\beta+1)\text{arccos}\sqrt{z}-\frac{\beta \pi}{2}-\frac{\pi}{4}) \\
                               -i\text{cos}((2n+\alpha+\beta-1)\text{arccos}\sqrt{z}-\frac{\beta \pi}{2}-\frac{\pi}{4})
                              \end{pmatrix}.
        \end{aligned}
       \end{equation*}
       Apparently, the fact is $R(z)=I+\mathcal{O}(n^{-1})$. The error terms hold uniformly for $x$ in compact subsets of $(0,1)$. Then, Theorem 2 is given by some calculations.

  \par As we have known in Section 2, $\Psi(\xi,\zeta)$ can be approximated by a Bessel model RH problem as $\zeta \to 0^+$, while it convert to a Airy function as $\zeta \to \infty$. For expressing the strong asymptotic expansion near the endpoint 1, we use the fact that the local parametrix can be approximated by the Bessel function  as $\zeta=2n^2t \to 0^+$.
  \\\textbf{Proof of Theorem 3.} After taking a series of inverse transformation of $Y(z) \to T(z) \to S(z) \to R(z)$ with the formulas (\ref{T}),\;(\ref{S}),\;(\ref{pinfty}),and (\ref{R}) as well as (\ref{p1}), we find,

  \begin{small}
    \begin{equation}\label{y11y225}
      \begin{aligned}
        \begin{pmatrix}
            Y_{11}(z)\\
            Y_{22}(z)
        \end{pmatrix} & =  4^{-n\sigma_3} R(z) P_-^{(1)}(z) \begin{pmatrix}
                         1&0\\
                         -w(z)^{-1}\varphi_-(z)^{-2n}  & 1
                        \end{pmatrix}\varphi_-(z)^{n\sigma_3}e^{\frac{t}{2z(1-z)}\sigma_3}\begin{pmatrix}
                          1 \\
                          0
                        \end{pmatrix}\\
                     & =   4^{-n\sigma_3} R(z) E_1(z) \Psi_-(\xi, \zeta) W_-(z)^{-\sigma_3}\varphi_-(z)^{-n\sigma_3}\\
                     & \begin{pmatrix}
                        1 & 0 \\
                        -w(z)^{-1}\varphi_-(z)^{-2n}  & 1
                       \end{pmatrix}\varphi_-(z)^{n\sigma_3}e^{\frac{t}{2z(1-z)}\sigma_3}\begin{pmatrix}
                         1 \\
                         0
                       \end{pmatrix}\\
                     & =   4^{-n\sigma_3} R(z) E_1(z) \Psi_-(\xi, \zeta)\begin{pmatrix}
                         e^{\frac{\beta \pi i}{2}}\\
                         -e^{-\frac{\beta \pi i}{2}}
                     \end{pmatrix} (w_{p_J2}(z,t))^{-\frac{1}{2}}
      \end{aligned}
     \end{equation}
   \end{small}
   where $\xi = n^2f_1(z)$ and $\zeta=2n^2t$.
   Then using the facts that
    \begin{equation}
        \begin{aligned}
            \varphi_-(z)=e^{-2i\text{arccos}\sqrt{z}} \text{and} \sqrt{n^2f_1(z)}_-= -2ni \text{arccos}\sqrt{z},
        \end{aligned}
    \end{equation}
    and the conclusion that  $I_{\theta}(z) e^{\frac{\theta \pi i}{2}}= J_\theta(ze^{\frac{\pi i }{2}})$ with $\text{arg} z \in (-\pi, \frac{\pi}{2}]$, see \cite{Olver2010NIST}, we get
  \begin{equation}\label{y11y224}
    \begin{aligned}
       \begin{pmatrix}
          Y_{11}(z)\\
           Y_{22}(z)
       \end{pmatrix} & =  4^{-n\sigma_3} R(z) E_1(z) R_{0}(\zeta) F_(\xi,\zeta) (w_{p_J2}(z,t))^{-\frac{1}{2}}             \begin{pmatrix}
        e^{\frac{\beta\pi i}{2}}\\
        -e^{-\frac{\beta \pi i}{2}}
        \end{pmatrix}  \\
                  & =   4^{-n\sigma_3} R(z) E_1(z) R_{0}(\zeta) \sqrt{\pi} e^{\frac{\zeta}{\xi_-}} (w_{p_J2}(z,t))^{-\frac{1}{2}}  \begin{pmatrix}
                   e^{\frac{\beta \pi i}{2}}I_{\beta}(\sqrt{\xi})_-\\
                   i  e^{\frac{\beta \pi i}{2}} \sqrt{\xi}_-  I'_\beta(\sqrt{\xi})_-
                    \end{pmatrix}  \\
                  & =    4^{-n\sigma_3} R(z) E_1(z) R_{0}(\zeta)\\
                  &\qquad \qquad \begin{pmatrix}
                    J_\beta(2n\text{arccos}\sqrt{z})\\
                    i2n\text{arccos}\sqrt{z} J'_\beta(2n\text{arccos}\sqrt{z})
                 \end{pmatrix} \sqrt{\pi}  e^{\frac{\zeta}{\xi_-}} (w_{p_J2}(z,t))^{-\frac{1}{2}}
    \end{aligned}
  \end{equation}
  where $J_\beta$ is the Bessel function of order $\beta$.
  To study the asymptotic behavior of $E_1(z)$ in (\ref{E1}), we compute the asymptotic behavior of $P^{(\infty)}$ from the upper plane. From the expression of $P^{(\infty)}$ in (\ref{pinfty}), we find
  \begin{equation}\label{Pi+}
    \begin{aligned}
        P_+^{(\infty)}(z)=\frac{e^{-\frac{\pi i}{4}}}{2x^{\frac{1}{4}}(1-x)^{\frac{1}{4}}} 2^{-(\alpha+\beta)\sigma_3} \begin{pmatrix}
            e^{i \theta_1} & i e^{-i \theta_1}  \\
            -ie^{i \theta_2} & e^{-i \theta_2}
           \end{pmatrix} e^{-\frac{\beta \pi i}{2}\sigma_3} x^{-\frac{\alpha}{2}\sigma_3} (1-x)^{-\frac{\beta}{2}\sigma_3}
    \end{aligned}
  \end{equation}
  where $\theta_1 =(\alpha+\beta+1)\text{arccos}\sqrt{z}$ and  $\theta_2 =(\alpha+\beta-1)\text{arccos}\sqrt{z}$.
  It is readily found that
  \begin{equation}
    \begin{aligned}
              (n^2f_1(z))_+^{\frac{1}{4}}=e^{\frac{\pi i }{4}\sigma_3} (2n\text{arccos}\sqrt{z})^{\frac{1}{2}\sigma_3},
    \end{aligned}
   \end{equation}
   Combining $E_1(z)$ in (\ref{E1}) with $W(z)$ in (\ref{W}), we obtain
   \begin{equation}\label{E11}
    \begin{aligned}
         E_1(z)= \frac{e^{-\frac{\pi i }{4}}}{\sqrt{2} x^{\frac{1}{4}}(1-x)^{\frac{1}{4}}} 2^{-(\alpha+\beta)\sigma_3}\begin{pmatrix}
                  \text{cos}\theta_1 & \text{sin}\theta_1\\
                  -i\text{cos}\theta_2& -i\text{sin}\theta_2
             \end{pmatrix}e^{\frac{\pi i}{4}\sigma_3}(2n\text{arccos}\sqrt{z})^{\frac{1}{2}\sigma_3}.
    \end{aligned}
   \end{equation}

   By substituting the formulas (\ref{E11}) into (\ref{y11y224}) and considering the result $Y_{11}(z)=\pi_n(z)$, $R(z)=I+\mathcal{O}(n^{-1})$, and $R_{0}(\zeta)=I+\mathcal{O}(\zeta)$, we prove the conclusion of Theorem 3, where the error terms is uniformly for $z\in \mathbb{D}_1$, as $\zeta = 2n^2t \to 0 $, and $\mathbb{D}_1=\{ x|2n\text{arccos}\sqrt{x}\in (0,\infty), n\to \infty \} $.
\par In order to obtain the strong asymptotic expansion near the endpoint 1, we take advantage of the fact $\Psi(\xi,\zeta)$ can be approximated by the Airy function when $\zeta=2n^2t \to \infty$.
  \\\textbf{Proof of Theorem 4.} It is convenient to start with (\ref{y11y225}), and we list it here
   \begin{equation*}
     \begin{aligned}
      \begin{pmatrix}
          Y_{11}(z)\\
          Y_{22}(z)
      \end{pmatrix}   & =   4^{-n\sigma_3} R(z) E_1(z) \Psi_-(\xi, \zeta)\begin{pmatrix}
        e^{\frac{\beta \pi i}{2}}\\
        -e^{-\frac{\beta \pi i}{2}}
     \end{pmatrix} (w_{p_J2}(z,t))^{-\frac{1}{2}}
     \end{aligned}
   \end{equation*}
    where $E_1(z)$ is in (\ref{E1}), $\xi=n^2 f_1(z)$, and $\zeta= 2n^2t$.
    \par Combining (\ref{H}), (\ref{N}), (\ref{U1}) and (\ref{R01}), we get
    \begin{equation}\label{Psi}
        \begin{aligned}
             \Psi(\zeta^{\frac{2}{3}}\lambda, \zeta)=\zeta^{-\frac{1}{6}\sigma_3}R_{01}(\lambda)E_{01}(\lambda)\Psi_{A}(\zeta^{\frac{2}{9}}f_{01}(\lambda))e^{\pm \frac{\beta \pi i}{2}\sigma_3}
        \end{aligned}
    \end{equation}
    where $\pm \text{Im} \lambda >0$, $\lambda \in U(-1,r)$, $\Psi_{A}$ is the Airy kernel. By integrating equations (\ref{U}), (\ref{R01p}) and (\ref{E01}), we acquire
    \begin{small}
    \begin{equation}\label{E01plus}
        \begin{aligned}
            E_{01}(\lambda)=(\lambda+1)^{-\frac{1}{4}\sigma_3} \frac{I+i\sigma_1}{\sqrt{2}} (\frac{\sqrt{\lambda+1}+1}{\sqrt{\lambda+1}-1})^{-\frac{\alpha}{2}\sigma_3} e^{\mp \frac{\alpha \pi i}{2}\sigma_3} \frac{I+i\sigma_1}{\sqrt{2}}^{-1} (\lambda^{\frac{2}{9}}f_{01}(\lambda))^{\frac{1}{4}\sigma_3}
        \end{aligned}
    \end{equation}
    \end{small}
    for $\pm \text{Im} \lambda>0$, $f_{01}(\lambda)=(-\frac{3}{2}\theta(\lambda))^{\frac{2}{3}}$, $\theta(\lambda)=\frac{(\lambda+1)^{\frac{3}{2}}}{\lambda}$.
    \par Then, after rescaling, we obtain,
    \begin{equation}\label{y11y226}
        \begin{aligned}
            \begin{pmatrix}
                Y_{11}(z)\\
                Y_{22}(z)
            \end{pmatrix}   & =   4^{-n\sigma_3} R(z) E_1(z) \zeta^{-\frac{1}{6}\sigma_3} R_{01}(\zeta) E_{01}(\lambda) \\
                            & \qquad \qquad \Psi_{A-}(\zeta^{\frac{2}{9}}f_{01}(\lambda)) \begin{pmatrix}
                                    1\\
                                   -1
                                  \end{pmatrix}(w_{p_J2}(z,t))^{-\frac{1}{2}}\\
                             & =   4^{-n\sigma_3} R(z) E_1(z) \zeta^{-\frac{1}{6}\sigma_3} R_{01}(\zeta) E_{01}(\lambda) \\
                             & \qquad \qquad \sqrt{2 \pi} \begin{pmatrix}
                                   Ai_-(\zeta^{\frac{2}{9}}f_{01}(\lambda))\\
                                   - i Ai'_-(\zeta^{\frac{2}{9}}f_{01}(\lambda))
                                  \end{pmatrix} (w_{p_J2}(z,t))^{-\frac{1}{2}}
        \end{aligned}
    \end{equation}
    By some simple calculations, it is easily to verify that
    \begin{equation}\label{lam12}
        \begin{aligned}
            (\lambda+1)_-^{-\frac{1}{4}\sigma_3}=(|\lambda|-1)^{-\frac{1}{4}\sigma_3}e^{\frac{\pi i}{4}\sigma_3},
        \end{aligned}
    \end{equation}
    \begin{equation}\label{lam1lam12}
        \begin{aligned}
            (\frac{\sqrt{\lambda+1}+1}{\sqrt{\lambda+1}-1})_-=\text{exp} (-i\;\text{arccos}\frac{1}{\sqrt{|\lambda|}}+i\pi),
        \end{aligned}
    \end{equation}
    \begin{equation}\label{f012}
        \begin{aligned}
            (f_{01}(\lambda))^{-\frac{3}{2}}=\frac{2}{3} (|\lambda|-1)^{-\frac{3}{2}} |\lambda| e^{\frac{3 \pi i}{2}}.
        \end{aligned}
    \end{equation}
    Substituting (\ref{lam12}), (\ref{lam1lam12}) and (\ref{f012}), we derive
    \begin{equation}\label{zeta1}
        \begin{aligned}
            \zeta^{-\frac{1}{6}\sigma_3}E_{01}(\lambda)=\begin{pmatrix}
                d_{11}& i d_{12}\\
                i d_{21} & d_{22}
            \end{pmatrix}
        \end{aligned}
    \end{equation}
    where $d_{11}$, $d_{12}$, $d_{21}$, and $d_{22}$ are defined as (\ref{d11}), (\ref{d12}), (\ref{d21}) and (\ref{d22}) respectively.
    With the help of the formula (\ref{R01pp}), it is readily to show, as $\zeta \to \infty$
    \begin{equation}\label{zeta2}
        \begin{aligned}
            \zeta^{-\frac{1}{6}\sigma_3}R_{01} \zeta^{\frac{1}{6}\sigma_3} = (I+\mathcal{O}(\zeta^{-\frac{1}{3}})\begin{pmatrix}
           1& 0   \\
           -\frac{7}{72(|\lambda|-1)}i & 1
            \end{pmatrix}.
        \end{aligned}
    \end{equation}
    \par Inserting the above results into (\ref{y11y226}) as well as the facts that $R_{01}(\zeta)=I+\mathcal{O}(\zeta^{-\frac{1}{3}})$,which is uniformly for $z \in \mathbb{D}_2$, as $\zeta = 2n^2 t \to \infty$ and $R(z)=I+\mathcal{O}(n^{-1})$ which are uniformly for $z \in (1-r,1)$, we follow Theorem 4.

  \par To consider the strong asymptotic expansions of the monic polynomial, we recall some results in \cite{ChenThe} that $\bar{\Psi}(\bar{\xi},\bar{\zeta})$ can be approximated by a modified Bessel model RH problem as $\bar{\zeta} \to 0$ and then prove the strong asymptotic expansion of monic OPS for $z \in (0,\delta)$, where $\delta$ is a small and positive constant.
   \\\textbf{Proof of Theorem 5.} As we say above, $\pi_n(z)$ is a polynomial function and we just consider the limit from the upper plane. Then, using the equations (\ref{T}),\;(\ref{S}),\;(\ref{pinfty}),\;and (\ref{R}), it is readily to show
   \begin{small}
   \begin{equation}\label{y11y222}
    \begin{aligned}
           \begin{pmatrix}
               Y_{11}(z)\\
               Y_{22}(z)
           \end{pmatrix}  &= Y_+(z)\begin{pmatrix}
                             1\\
                             0
                          \end{pmatrix} \\
                          & =  4^{-n\sigma_3} R(z) P_+^{(0)}(z) \begin{pmatrix}
                                                                   1 & 0   \\
                                                                   w(z)^{-1}\varphi_+(z)^{-2n} &1
                                                               \end{pmatrix} e^{\frac{t}{2z(1-z)}\sigma_3}\varphi_+(z)^{n\sigma_3}\\
                          & = 4^{-n\sigma_3} R(z)E_0(z) \bar{\Psi}_-(\bar{\xi}, \bar{\zeta}) e^{\frac{(\alpha-1)\pi i}{2}\sigma_3} \begin{pmatrix}
                                1& 0   \\
                                1&1
                           \end{pmatrix} (w_{p_J2}(z,t))^{-\frac{1}{2}\sigma_3}  \begin{pmatrix} 1\\0 \end{pmatrix},
    \end{aligned}
   \end{equation}
   \end{small}
    where $w(z)=z^\alpha (1-z)^\beta$, $\bar{\xi} = n^2f_0(z)$, $\bar{\zeta} = 2n^2 t$, and $E_0$ is in (\ref{E0}). Besides, the boundary value on the positive side $P_+^{(0)}$ corresponds to the value $\bar{\Psi}_-$ on the negative side, as can be seen from the correspondence $\bar{\xi} \approx -4n^2z$ when $z \to 0$.
    According to the results of Chen has proven, see (86),(87) and (89) in \cite{ChenThe}, then
   \begin{small}
    \begin{equation}\label{y11y221}
        \begin{aligned}
           \begin{pmatrix}
              Y_{11}(z)\\
               Y_{22}(z)
           \end{pmatrix} & =  4^{-n\sigma_3} R(z) E_0(z) R_{0S}(\bar{\zeta}) G_(\bar{\xi},\bar{\zeta}) (w_{p_J2}(z,t))^{-\frac{1}{2}}             \begin{pmatrix}
            e^{\frac{(\alpha-1)\pi i}{2}}\\
            e^{-\frac{(\alpha-1)\pi i}{2}}
            \end{pmatrix}  \\
                      & =   4^{-n\sigma_3} R(z) E_0(z) R_{0S}(\bar{\zeta})(-1)^n \sqrt{\pi} e^{\frac{\bar{\zeta}}{\bar{\xi}_-}} (w_{p_J2}(z,t))^{-\frac{1}{2}}  \begin{pmatrix}
                       -iJ_\alpha(\sqrt{|\bar{\xi}|})\\
                       \sqrt{|\bar{\xi}|} J'_\alpha(\sqrt{|\bar{\xi}|})
                        \end{pmatrix}
        \end{aligned}
    \end{equation}
    \end{small}
    accompanied with the facts that
    \begin{equation}
        \begin{aligned}
          \bar{\xi}_-=n^2 f_{0-}(z)= e^{-\pi i}n^2 (\pi - 2\text{arccos}\sqrt{z})^2
        \end{aligned}
    \end{equation}
    and
    \begin{equation}
        \begin{aligned}
            \sqrt{|\bar{\xi}|}=n(\pi - 2\text{arccos}\sqrt{z}).
        \end{aligned}
    \end{equation}
    Also, we give a more precise expression for $E_0(z)$ by the formula (\ref{E0}), (\ref{pinfty}), (\ref{V}) and (\ref{Pi+})
    \begin{equation}\label{E0+}
        \begin{aligned}
            E_0(z)= &\frac{(-1)^n e^{-\frac{\pi i}{4}}}{\sqrt{2}x^{\frac{1}{4}}(1-x)^{\frac{1}{4}}} 2^{-(\alpha+\beta)\sigma_3}\begin{pmatrix}
                    -\text{sin}\theta_3 & \text{cos}\theta_3 \\
                    i\text{sin}\theta_4 & -i \text{cos}\theta_4
                   \end{pmatrix}\\
                   & e^{-\frac{\pi i}{4}\sigma_3} (n(\pi  - 2\text{arccos}\sqrt{z}))^{\frac{1}{2}\sigma_3}
        \end{aligned}
    \end{equation}
   \par Bringing the results of (\ref{Pi+}) and (\ref{E0+}) back to (\ref{y11y221}), we obtain the result of Theorem 5. What's more, $R(z)=I+\mathcal{O}(n^{-1})$ uniformly for $x \in (0, \delta)$ with $0<\delta \ll \frac{1}{2}$ and $R_{0S}(\bar{\zeta})=I+\mathcal{O}(\bar{\zeta})$ uniformly for $z \in \mathbb{D}_3$ as $\bar{\zeta}= 2n^2 t \to 0$, where $\mathbb{D}_3=\{ x|n(\pi-2\text{arccos}\sqrt{x})\in (0,\infty), n\to \infty\}$.

   We continue to consider the strong asymptotic expansion for $z \in (0, \varepsilon)$, by using the theorem in \cite{ChenThe} that the model RH problem for $\bar{\Psi}(\bar{\xi},\bar{\zeta})$ can be approximated by the Airy model when $\bar{\zeta} = 2n^2 t \to \infty$.
   \\ \textbf{Proof of Theorem 6.} For your convenience, this application start with (\ref{y11y222}), that is
   \begin{small}
   \begin{equation}
    \begin{aligned}
           \begin{pmatrix}
               Y_{11}(z)\\
               Y_{22}(z)
           \end{pmatrix}
                           = 4^{-n\sigma_3} R(z)E_0(z) \bar{\Psi}_-(\bar{\xi},\bar{\zeta}) e^{\frac{(\alpha-1)\pi i}{2}\sigma_3} \begin{pmatrix}
                                1& 0   \\
                                1&1
                           \end{pmatrix} (w_{p_J2}(z,t))^{-\frac{1}{2}\sigma_3}  \begin{pmatrix} 1\\0 \end{pmatrix},
    \end{aligned}
   \end{equation}
   \end{small}
    where $E_0(z)$ is in (\ref{E0}), $\bar{\xi} = n^2 f_0(z)$ and $\bar{\zeta}= 2n^2 t$.
    \par Then, we use the same way as we deal with the situation for $z \in (1-r,1)$. Thus, we rescale $\bar{\xi}$ by $\bar{\zeta}^{\frac{2}{3}}\bar{\lambda}$ with the help of the equation (92), (93), (95) and (96) in \cite{XuCritical}(To avoid the symbol confusion, we add a superscript in these equation), and we obtain
    \begin{equation}\label{y11y223}
        \begin{aligned}
            \begin{pmatrix}
                Y_{11}(z)\\
                Y_{22}(z)
            \end{pmatrix} & = 4^{-n\sigma_3} R(z) E_0(z) \bar{\zeta}^{-\frac{1}{6}\sigma_3}\bar{R}_{01}\bar{E}_{01}(\bar{\lambda})\\
                          & \qquad \qquad \bar{\Psi}_{A-}(\bar{\zeta}^{\frac{2}{9}}f_{01}(\bar{\lambda}))(-1)^{n\sigma_3}\begin{pmatrix}
                             -i\\
                              i
                             \end{pmatrix}(w_{p_J2}(z,t))^{-\frac{1}{2}}\\
                          & = 4^{-n\sigma_3} R(z) E_0(z) \bar{\zeta}^{-\frac{1}{6}\sigma_3}\bar{R}_{01}\bar{E}_{01}(\bar{\lambda})\\
                        & \qquad \qquad \begin{pmatrix}
                             -i Ai_-(\bar{\zeta}^{\frac{2}{9}}f_{01}(\bar{\lambda}))\\
                             - Ai'_-(\bar{\zeta}^{\frac{2}{9}}f_{01}(\bar{\lambda}))
                            \end{pmatrix}(-1)^n \sqrt{2 \pi}(w_{p_J2}(z,t))^{-\frac{1}{2}}
        \end{aligned}
    \end{equation}

    \par Finally, we can easily find the similar properties in (\ref{zeta1}) and (\ref{zeta2}) for $\bar{E}_{01}$ and $\bar{R}_{01}$. By substituting these formulas into (\ref{y11y223}), we have reached the desired conclusion. Besides, $R(z)=I+\mathcal{O}(n^{-1})$ holds uniformly for $z \in (0,\epsilon)$ while $\bar{R}_{01}= I + \mathcal{O}(\bar{\zeta}^{-\frac{1}{3}})$ is uniformly for $z \in \mathbb{D}_4$, see (\ref{D4}), as $\bar{\zeta}= 2n^2t \to \infty$.
    \par The next step of this paper is to proof the limit behavior of the kernel $K_n(x,y)$ in the bulk of the spectrum.
   \\ \textbf{Proof of Theorem 7.}  For convenience, we start with (\ref{kn}) and rewrite it here:
    \begin{equation}
        \begin{aligned}
            K_n(x,y;t)=\gamma_{n-1}^2\sqrt{w_{p_J2}(x,t)}\sqrt{w_{p_J2}(y,t)} \frac{\pi_n(x)\pi_{n-1}(y)-\pi_n(y)\pi_{n-1}(x)}{(x-y)}.
        \end{aligned}
    \end{equation}
    \par Then,  by the facts that $Y_{11}(z)=\pi_n(z)$ and $\pi_n(z)$ is analytic in the whole complex plane, one finds
    \begin{equation}\label{knY}
        \begin{aligned}
            K_n(x,y;t)& =\frac{\sqrt{w_{p_J2}(x,t)}\sqrt{w_{p_J2}(y,t)}}{2\pi i (x-y)}(Y_+^{-1}(y)Y_+(x))_{21}\\
                      & = \frac{\sqrt{w_{p_J2}(x,t)}\sqrt{w_{p_J2}(y,t)}}{2\pi i (x-y)}\begin{pmatrix}
                          0 & 1
                      \end{pmatrix}Y_+^{-1}(y)Y_+(x)\begin{pmatrix}
                            1\\
                            0
                      \end{pmatrix}.
        \end{aligned}
    \end{equation}
    \par In order to express $Y(z)$, we take the inverse procedure of $Y \to T \to S \to R$ in (\ref{T}),\;(\ref{S})\;and (\ref{R}) and get
    \begin{equation}\label{yplus}
        \begin{aligned}
            Y_+(x) & = 4^{-n\sigma_3} R(x) P^{(\infty)}_+(x) \begin{pmatrix}
                      1& 0\\
                      w(x)^{-1}\varphi_+(x)^{-2n} &1
                     \end{pmatrix} e^{\frac{t}{2x(1-x)_+}\sigma_3}\varphi_+(x)^{n\sigma_3}\\
                  & = 4^{-n\sigma_3} R(x) 2^{-(\alpha+\beta)\sigma_3}M^{-1} a(x)^{-\sigma_3} M D_+(x)^{\sigma_3}\\
                  & \begin{pmatrix}
                    1& 0\\
                    w(x)^{-1}\varphi_+(x)^{-2n} &1
                   \end{pmatrix} e^{\frac{t}{2x(1-x)_+}\sigma_3}\varphi_+(x)^{n\sigma_3}
        \end{aligned}
    \end{equation}
    where $M = \frac{I+i\sigma_1}{\sqrt{2}}$, $a(z)=(\frac{z-1}{z})^{\frac{1}{4}} $, $ D_+(x)=\varphi_+(x)^{\frac{\alpha+\beta}{2}}x^{-\frac{\alpha}{2}}(x-1)_+^{-\frac{\beta}{2}}$ and $R(z)=I+\mathcal{O}(n^{-\frac{2}{3}})$. By the facts that,
   \begin{equation*}
    \begin{aligned}
              \varphi_+(x)= e^{2i \text{arccos}\sqrt{x}},
    \end{aligned}
   \end{equation*}
    we find,
    \begin{equation}\label{yw}
        \begin{aligned}
            Y_+(x) (w_{p_J2}(x,t))^{\frac{1}{2}}    \begin{pmatrix}
                1\\
                0
            \end{pmatrix}& = 4^{-n\sigma_3} R(x) 2^{-(\alpha+\beta)\sigma_3} M^{-1} a(x)^{-\sigma_3} M \\
                      &\begin{pmatrix}
                       \text{exp}(-\frac{\beta \pi i}{2}+i(\alpha+\beta+2n)\text{arccos}\sqrt{x})\\
                      \text{exp}(\frac{\beta \pi i}{2}-i(\alpha+\beta+2n)\text{arccos}\sqrt{x})
                     \end{pmatrix}
        \end{aligned}
    \end{equation}

    \par Moreover, it can be easily proved that
   \begin{equation}\label{aa}
    \begin{aligned}
         M^{-1}a(y)^{\sigma_3}a(x)^{-\sigma_3}M= I +\mathcal{O}(x-y),
    \end{aligned}
   \end{equation}
   where the error terms  holds uniformly for any $x,\;y$ in the compact subsets of $(0,1)$.
  Inserting (\ref{yw}) and (\ref{aa}) into (\ref{yplus}) , one finds
    \begin{equation}\label{knxy}
        \begin{aligned}
            K_n(x,y;t)= \frac{\text{sin}[(2n+\alpha+\beta)(\text{arccos}\sqrt{y}-\text{arccos}\sqrt{x})]}{\pi (x-y)}+\mathcal{O}(1).
        \end{aligned}
    \end{equation}
    If $x \to y$, then
    \begin{equation*}
        \begin{aligned}
            \frac{1}{n}K_n(y,y) =\frac{1}{\pi \sqrt{y(1-y)}}+\mathcal{O}(\frac{1}{n}) \qquad n \to \infty,
        \end{aligned}
    \end{equation*}
   for any $y \in (0,1)$ with the uniform error terms.
   \par Let $a \in (0,1)$, $u,v\in \mathbb{R}$, and $\rho=(\pi(\sqrt{x(1-x)}))^{-1}$. Then, $x$ is rescaled as $a + \frac{u}{n\rho(a)}$ while $y$ is rescaled as $a + \frac{v}{n\rho(a)}$ in (\ref{knxy}). After some calculations, we obtain
   \begin{small}
   \begin{equation}\label{2nab}
    \begin{aligned}
        (2n+\alpha+\beta)(\text{arccos}\sqrt{a + \frac{v}{n\rho(a)}}-\text{arccos}\sqrt{a + \frac{u}{n\rho(a)}} )= \pi(u-v)+\mathcal{O}(\frac{1}{n}),
    \end{aligned}
   \end{equation}
   \end{small}
   uniformly for $a$ in any compact subsets of $(0,1)$. Thus, combining (\ref{knxy}) with (\ref{2nab}), the result can be proved.

   For the sake of expressing the limiting kernel at the right edge $1$, we define the $\psi$-functions
   \begin{equation}\label{psi121}
    \begin{aligned}
        \begin{pmatrix}
            \psi_1(\xi,\zeta)\\
            \psi_2(\xi,\zeta)
        \end{pmatrix}= \Psi_+(\xi, \zeta)\begin{pmatrix}
            e^{-\frac{\beta \pi i}{2}}\\
            e^{\frac{\beta \pi i}{2}}
        \end{pmatrix}, \quad  \xi<0.
    \end{aligned}
    \end{equation}
   The functions are determined via a model RH problem $\Psi(\xi,\zeta)$ related to a special solution of a third-order nonlinear differential equation which can be reduced to a certain Painlev$\acute{e}$ \uppercase\expandafter{\romannumeral3} equation, see \cite{XuCritical}.
   We also give the expression of  the $\Psi$-kernel as
    \begin{equation}\label{kp}
        \begin{aligned}
            K_\Psi (u,v,\zeta)= \frac{\psi_1(-u,\zeta) \psi_2(-v,\zeta)-\psi_1(-v,\zeta)\psi_2(-u,\zeta)}{2\pi i (u-v)}
        \end{aligned}
    \end{equation}
    for $u,v,\zeta \in (0,\infty) $.
    It is convenient to begin with (\ref{knY}).
    \begin{equation}\label{kn1}
    \begin{aligned}
        K_n(x,y;t)
                   = \frac{\sqrt{w_{p_J2}(x,t)}\sqrt{w_{p_J2}(y,t)}}{2\pi i (x-y)}\begin{pmatrix}
                      0 & 1
                  \end{pmatrix}Y_+^{-1}(y)Y_+(x)\begin{pmatrix}
                        1\\
                        0
                  \end{pmatrix}.
    \end{aligned}
    \end{equation}
   \par  Tracing back the transformation $R \to S \to T \to Y$ with the formulas (\ref{T}),\;(\ref{S}),\;(\ref{pinfty}),and (\ref{R}), one finds
    \begin{equation}
        \begin{aligned}
            Y_+(z) & = 4^{-n\sigma_3} R(z) P^{(1)}_+(z) \begin{pmatrix}
                 1 & 0\\
                 w(z)^{-1}\varphi_+(z)^{-2n} &1
            \end{pmatrix} e^{\frac{t}{2z(1-z)}\sigma_3} \varphi_+(z)^{n\sigma_3}\\
                 & = 4^{-n\sigma_3} R(z) E_1(z) \Psi_+(\xi,\zeta) e^{-\frac{\beta \pi i}{2}\sigma_3}\begin{pmatrix}
                     1&0\\
                     1&1
                 \end{pmatrix}w_{p_J2}(z,t)^{-\frac{1}{2}}
        \end{aligned}
    \end{equation}
    Accordingly, we obtain
    \begin{equation}\label{yw1}
        \begin{aligned}
            Y_+(x)w_{p_J2}(x,t)^{\frac{1}{2}}\begin{pmatrix}
                1\\
                0
            \end{pmatrix} & =4^{-n\sigma_3} R(x) E_1(x) \Psi_+(\xi,\zeta) e^{-\frac{\beta \pi i}{2}\sigma_3}\begin{pmatrix}
                             1\\
                             1
                            \end{pmatrix}\\
                         & = 4^{-n\sigma_3} R(x) E_1(x) \Psi_+(\xi,\zeta)\begin{pmatrix}
                            e^{-\frac{\beta \pi i}{2}}\\
                            e^{\frac{\beta \pi i}{2}}
                         \end{pmatrix}.
        \end{aligned}
    \end{equation}
    \par Inserting (\ref{yw1}) into (\ref{kn1}), it is easy to achieve
    \begin{small}
        \begin{equation}\label{kn3}
            \begin{aligned}
                K_n(x,y)=\frac{(-\psi_2(g_n(y)),\psi_1(g_n(y)))E_1(y)^{-1}R(y)^{-1}R(x)E_1(x)(\psi_1(g_n(x)),\psi_2(g_n(x)))^{T}}{2\pi i (x-y)}
            \end{aligned}
        \end{equation}
       \end{small}
    where $g_n(z)=n^2f_1(z) $.
    \par Let $x=1-\frac{u}{4n^2}$ and $y=1-\frac{v}{4n^2}$ with $u,v = \mathcal{O}(1)$. Besides, due to
       $$g_n(z)=n^2f_1(z)=n^2[4(z-1)+\mathcal{O}(z-1)^2], \quad z \to 1,$$
    we have
    \begin{equation}\label{gn}
        \begin{aligned}
            g_n(x)=-u(I+\mathcal{O}(\frac{1}{n^2})),\quad g_n(y)=-v(I+\mathcal{O}(\frac{1}{n^2}))
        \end{aligned}
       \end{equation}
    \par Since $E_1(z)$ is a matrix function analytic in $U(1,r)$, it is directly to show
    \begin{equation}
        \begin{aligned}
            E_1(y)^{-1}E_1(x) & = I +E_1(y)^{-1}(E_1(x)-E_1(y))\\
                              &= I+\mathcal{O}(x-y)\\
                              &= I+(u-v)\mathcal{O}(n^{-2})
        \end{aligned}
     \end{equation}
   for bounded $u,v$. Similarly, the analyticity of $R(z)$ for $z \in U(1,r)$ implies that
   \begin{equation}
    \begin{aligned}
        R^{-1}(y)R(x)=I+\mathcal{O}(x-y)=I+(u-v)\mathcal{O}(n^{-2}),
    \end{aligned}
   \end{equation}
   here again with uniform error terms. Therefore, by the above formulas, we obtain
   \begin{small}
    \begin{equation}
     \begin{aligned}
         K_n(x,y)=\frac{(-\psi_2(g_n(y)),\psi_1(g_n(y)))(I+\mathcal{O}(x-y))(\psi_1(g_n(x)),\psi_2(g_n(x)))^T}{2\pi i (x-y)}
     \end{aligned}
    \end{equation}
    \end{small}
   for $x,y \in (1-r, 1)$ where $r >0$.
   Meanwhile, the equation (\ref{gn}) implies that
   \begin{equation}
    \begin{aligned}
        \psi_k(g_n(x),\zeta)=\psi_k(-u,\zeta)+\mathcal{O}(n^{-2})\; \text{and} \; \psi_k(g_n(y),\zeta)=\psi_k(-v,\zeta)+\mathcal{O}(n^{-2})
    \end{aligned}
    \end{equation}
    for $k=1,2$, where the uniform error terms hold uniformly for $u$ and $v$ lying in compact subsets of $(0, \infty)$.
   \par Hence, combining the above formula, we have
   \begin{small}
    \begin{equation}
        \begin{aligned}
        \frac{1}{4n^2}K_n(1-\frac{u}{4n^2},1-\frac{v}{4n^2})=\frac{\psi_1(-u,\zeta)\psi_2(-u,\zeta)-\psi_1(-u,\zeta)\psi_2(-v,\zeta)}{2 \pi i(u-v)}+\mathcal{O}(\frac{1}{n^2}),
        \end{aligned}
     \end{equation}
    \end{small}
    as $n \to \infty$.
    \par It is also direct to proof the $\Psi$-kernel is reduced to the Bessel kernel when $\zeta \to 0^+$. After we substitute (\ref{Fxz}) and (\ref{R_0}) into (\ref{psi121}), we find
    \begin{equation}
        \begin{aligned}
            \begin{pmatrix}
                \psi_1(\xi,\zeta)\\
                \psi_2(\xi,\zeta)
            \end{pmatrix} & =R_{0}(\xi) \pi^{\frac{1}{2}\sigma_3}\begin{pmatrix}
                I_\beta(\sqrt{|\xi|}e^{\frac{\pi i}{2}})\\
                i \pi \sqrt{|\xi|} I'_\beta(\sqrt{|\xi|}e^{\frac{\pi i}{2}})
                \end{pmatrix}\\
                        & = R_{0}(\xi) \pi^{\frac{1}{2}\sigma_3}\begin{pmatrix}
                J_\beta(\sqrt{|\xi|})\\
                i \pi \sqrt{|\xi|} J'_\beta(\sqrt{|\xi|})
                \end{pmatrix},
        \end{aligned}
    \end{equation}
    During this process, we have use the equation $e^{-\frac{1}{2}\beta \pi i}I_\beta(z)=J_\beta(z e^{-\frac{1}{2}\pi i}), \; \text{arg} z \in (0, \frac{\pi}{2}]$ (see (10.27.6) in \cite{Olver2010NIST}), and $X^T \begin{pmatrix} 0&1\\ -1&0   \end{pmatrix} X = \begin{pmatrix} 0&1\\ -1&0   \end{pmatrix}$ as det$X=1$, wher $X^T$ means the transpose of a matrix $X$. It follows that, as $\zeta \to 0^+$,
        \begin{equation}
            \begin{aligned}
                \frac{1}{4n^2}K_n(1-\frac{u}{4n^2},1-\frac{v}{4n^2}) & =\frac{\psi_1(-u,\zeta)\psi_2(-v,\zeta)-\psi_1(-v,\zeta)\psi_2(-u,\zeta)}{2\pi i (u-v)}+\mathcal{O}(\frac{1}{n^2})\\
                                                                     & = \mathbb{J}_\beta (u,v)+\mathcal{O}(\zeta)+\mathcal{O}(1/n^2).
            \end{aligned}
         \end{equation}
    where $\mathbb{J}_\beta$ is the Bessel kernel, and the error terms are uniform in compact subsets of $u,v \in (0, \infty)$.
    Now we consider the asymptotic behaviors for $\Psi$-kernel as $\zeta \to \infty$, which transfers to the Airy kernel. Taking a series of transformations $\Psi(\xi,\zeta) \to H(\lambda,\zeta) \to N(\lambda,\zeta) \to R_{01}(\lambda)$ and $R_{01}(\lambda)= I+\mathcal{O}(\zeta^{-\frac{1}{3}}\lambda^{-1}),\; \zeta \to \infty, \;\lambda \to \infty$, we obtain, as $\zeta \to \infty,\; \lambda \to \infty$
    \begin{equation}
       \begin{aligned}
        \Psi(\zeta^{\frac{2}{3}}\lambda, \zeta)=\zeta^{-\frac{1}{6}\sigma_3}(I+\mathcal{O}(\zeta^{-1/3}\lambda^{-1}))U(\lambda)e^{\frac{1}{3}\theta(\lambda)\sigma_3}.
       \end{aligned}
    \end{equation}
   With the help of the approximation of $U(\lambda)$ as $\lambda \to \infty$, it gives that, as $\zeta \to \infty,\; \lambda \to \infty$
   \begin{equation}
    \begin{aligned}
        \Psi(\zeta^{\frac{2}{3}}\lambda, \zeta) e^{-\zeta^{\frac{1}{3}}\sqrt{\lambda}\sigma_3}= (\zeta^{\frac{2}{3}}\lambda)^{-\frac{1}{4}\sigma_3}[I + \mathcal{O}(\zeta^{\frac{1}{3}}\lambda^{\frac{1}{2}})]M
    \end{aligned}
   \end{equation}
    \par For $|\lambda+1|<r$, we begin with the formula (\ref{Psi}) and substitute this formula into (\ref{psi121}),
    \begin{equation}
        \begin{aligned}
            \begin{pmatrix}
                \psi_1(\xi,\zeta)\\
                \psi_2(\xi,\zeta)
            \end{pmatrix} & =\zeta^{-\frac{1}{6}\sigma_3}R_{01}(\lambda)E_{01}(\lambda) \Psi_A(\zeta^{\frac{2}{9}}f_{01}(\lambda)) e^{\frac{1}{2}\beta \pi i \sigma_3}\begin{pmatrix}
                e^{-\frac{\beta \pi i}{2}}\\
                e^{\frac{\beta \pi i}{2}}
                \end{pmatrix}\\
                        & = \zeta^{-\frac{1}{6}\sigma_3}R_{01}(\lambda)E_{01}(\lambda) \sqrt{2 \pi }\begin{pmatrix}
                Ai(\zeta^{\frac{2}{9}}f_{01}(\lambda))\\
               - i Ai'(\zeta^{\frac{2}{9}}f_{01}(\lambda))
                \end{pmatrix},
        \end{aligned}
    \end{equation}
    Let $\lambda=-1+\frac{u}{m\zeta^{2/9}}$, then it is easy to verify that $\zeta^{2/9}f_{01}(\lambda)=u[1+\mathcal{O}(u\zeta^{-\frac{2}{9}})]$. Thus, the approximation of $K_\Psi$ can be given by (\ref{kpsinfity}), and $K_n$ followed by (\ref{knsinfity}) as $\zeta \to +\infty$.
    By now, we have finished the proof of case 4 in Theorem 7.

   \par Then we focus on the large-$n$ behavior of the kernel near the edge $x=0$. We also begin with the equation (\ref{kn1}).
   In this case, $Y(z)$ has been found in (\ref{y11y222}), and after a compilation of equations, we show, for $x \in (0,d)$
    \begin{equation}\label{ypp}
        \begin{aligned}
            Y_+(x) (w_{p_J2}(x,t))^{\frac{1}{2}}    \begin{pmatrix}
                1\\
                0
            \end{pmatrix}& = 4^{-n\sigma_3} R(x) E_0(x) \bar{\Psi}_-(f_n(x),\zeta) e^{\frac{(\alpha-1)\pi i }{2}\sigma_3}\begin{pmatrix}
                           1\\
                           1
                          \end{pmatrix} \\
                         & =   4^{-n\sigma_3} R(z) E_0(x) \bar{\Psi}_-(f_n(x),\zeta)\begin{pmatrix}
                            e^{\frac{(\alpha-1) \pi i}{2}}\\
                            e^{-\frac{(\alpha-1) \pi i}{2}}
                          \end{pmatrix}
        \end{aligned}
    \end{equation}
    where $f_n(x)=n^2f_0(x)$ and $\zeta = 2n^2t$.
    \par We rewrite the fact (4.4) in \cite{xu2015critical}, for $\bar{\xi} \in (-\infty,0)$,
    \begin{equation}\label{psi12}
        \begin{aligned}
             \begin{pmatrix}
                \bar{\psi}_1(\bar{\xi})\\
                \bar{\psi}_2(\bar{\xi})
                \end{pmatrix}= \begin{pmatrix}
                    \bar{\psi}_1(\bar{\xi},\zeta)\\
                    \bar{\psi}_2(\bar{\xi},\zeta)
                    \end{pmatrix}=\bar{\Psi}_-(\bar{\xi},\zeta)\begin{pmatrix}
                        e^{\frac{(\alpha-1) \pi i}{2}}\\
                        e^{-\frac{(\alpha-1) \pi i}{2}}
                        \end{pmatrix}
        \end{aligned}
    \end{equation}
   Substituting (\ref{ypp}) and (\ref{psi12}) into (\ref{kn1}), we get
   \begin{small}
    \begin{equation}\label{kn2}
        \begin{aligned}
            K_n(x,y)=\frac{(-\bar{\psi}_2(f_n(y)),\bar{\psi}_1(f_n(y)))E_0(y)^{-1}R(y)^{-1}R(x)E_0(x)(\bar{\psi}_1(f_n(x)),\bar{\psi}_2(f_n(x)))^{T}}{2\pi i (x-y)}
        \end{aligned}
    \end{equation}
   \end{small}
    Now specifying
   \begin{equation}
    \begin{aligned}
             x = \frac{u}{4n^2},\quad y=\frac{v}{4n^2},
    \end{aligned}
   \end{equation}
   where $u,v \in (0, \infty)$ and of the size $O(1)$. Also, we know
      $$f_n(z)=n^2 f_0(z)=-4n^2z(1+\mathcal{O}(z)),\; \text{for} z \to 0 $$
   Then it follows that
   \begin{equation}\label{fn}
    \begin{aligned}
        f_n(x)=-u(I+\mathcal{O}(\frac{1}{n^2})),\quad f_n(y)=-v(I+\mathcal{O}(\frac{1}{n^2}))
    \end{aligned}
   \end{equation}
   The matrix function $E_0(z)$ is analytic in $U(0,d)$, so
   \begin{equation}\label{IEE}
    \begin{aligned}
        E_0^{-1}(y)E_0(x) & =I +E_0^{-1}(y)(E_0(x)-E_0(y))\\
                          &= I+\mathcal{O}(x-y)\\
                          &= I+(u-v)\mathcal{O}(n^{-2})
    \end{aligned}
   \end{equation}
   for the above $u,v$. Similarly, the analyticity of $R(z)$ in $U(0,d)$ tells us
   \begin{equation}\label{IRR}
    \begin{aligned}
        R^{-1}(y)R(x)= I+(u-v)\mathcal{O}(n^{-2})
    \end{aligned}
   \end{equation}
   and also the error term holds uniformly for $u,v \in (0,\infty)$.
   \par Hence, substituting (\ref{IEE}) and (\ref{IRR}) into (\ref{kn2}), we obtain
   \begin{small}
   \begin{equation}
    \begin{aligned}
        K_n(x,y)=\frac{(-\bar{\psi}_2(f_n(y)),\bar{\psi}_1(f_n(y)))(I+\mathcal{O}(x-y))(\bar{\psi}_1(f_n(x)),\bar{\psi}_2(f_n(x)))^T}{2\pi i (x-y)}
    \end{aligned}
   \end{equation}
   \end{small}
    the error term is uniformly for $x,y \in (0,d)$, with $d>0$ being constant.
    Also the formula in (\ref{fn}) implies that
    \begin{equation}
        \begin{aligned}
            \bar{\psi}_k(f_n(x),\zeta)=\bar{\psi}_k(-u,\zeta)+\mathcal{O}(n^{-2})
        \end{aligned}
    \end{equation}
    for $k=1,2$, again with uniform error terms.
    \par Therefore, we have
    \begin{small}
    \begin{equation}
        \begin{aligned}
        \frac{1}{4n^2}K_n(\frac{u}{4n^2},\frac{v}{4n^2})=\bar{K}_{\bar{\Psi}}(u,v;\zeta)+\mathcal{O}(\frac{1}{n^2})
        \end{aligned}
     \end{equation}
    \end{small}
    for $n \to \infty$, where
    \begin{equation}
        \begin{aligned}
             \bar{K}_{\bar{\Psi}}(u,v;\zeta)=\frac{\bar{\psi}_1(-u,s)\bar{\psi}_2(-v,s)-\bar{\psi}_1(-v,s)\bar{\psi}_2(-u,s)}{2 \pi i(u-v)}
        \end{aligned}
    \end{equation}
    is the Painlevé type kernel. The error term $\mathcal{O}(n^{-2})$ in uniform for $u,v$ in any compact subsets of $(0,\infty)$ and $\zeta \in (0,\infty)$, thus finishing the proof of Theorem 7.
    \begin{remark}
        In the study of \cite{XuCritical} and \cite{xu2015critical}, we obtain the properties of $\bar{\Psi}$-kernel. To be specific, $\bar{\Psi}$-kernel is approximated by the Bessel kernel as $\zeta \to 0^+$
      \begin{equation*}
        \begin{aligned}
            \bar{K}_{\bar{\Psi}} (u,v,\zeta)=\mathbb{J}_\alpha(u,v)+\mathcal{O}(\zeta),
        \end{aligned}
      \end{equation*}
      uniformly for $u,v$ in compact subsets of $(0,\infty)$. $\mathbb{J}_\alpha$ is the Bessel kernel. Using this property, the Bessel kernel limit of $K_n$ has been given under some conditions. Similarly, the $\bar{\Psi}$-kernel can be approximated by the Airy kernel as $\zeta \to +\infty$ and the Airy kernel limit has been expressed with some assumptions using the same method as we deal with the $\Psi$-kernel. Thus, we can get the similar results for the $\bar{\Psi}$-kernel and $K_n$ with different parameters.\vspace{8mm}
    \end{remark}

\noindent\textbf{Acknowledgements}

This work is supported by  the National Science
Foundation of China (Grant No. 11671095,  51879045).


\end{document}